\documentclass[12pt]{article}

\usepackage{verbatim}
\usepackage{amssymb}
\usepackage{latexsym}
\usepackage{color}
\definecolor{green}{rgb}{0,1,0}
\usepackage{cancel}
\newcommand*\bulletlike{
\begin{picture}(0,0)
\unitlength1pt
\put(0,0){\circle*{2.5}}
\end{picture}}

\usepackage{xurl}
\usepackage{xcolor}
\definecolor{linkColor}{RGB}{16,16,16}
\usepackage{hyperref}
\hypersetup{colorlinks,allcolors=linkColor}

\hypersetup{
           breaklinks=true,   
           colorlinks=true,   
           pdfusetitle=true,  
        }
\begin{document}
\newfont{\blb}{msbm10 scaled\magstep1} 

\newtheorem{theo}{Theorem}[section]

\newtheorem{defi}[theo]{Definition}
\newtheorem{prop}[theo]{Proposition}
\newtheorem{lemm}[theo]{Lemma}
\newtheorem{coro}[theo]{Corollary}
\pagestyle{myheadings}
\date{}
\author{Layla Sorkatti \\
Gunnar Traustason \\
Department of Mathematical Sciences \\ 
University of Bath, UK}

\title{Nilpotent symplectic alternating algebras III}
\maketitle
\begin{abstract}
\mbox{}\\
In this paper we finish our classification of nilpotent symplectic alternating 
algebras of dimension 10 over any field $\mathbb{F}$.
\end{abstract}

\section{Introduction}
A symplectic alternating algebras (SAA) is a symplectic vector space L, whose associated alternating form is non-degenerate, that is furthermore equipped with a binary alternating product $\cdot : L \times L \mapsto L$ with the extra requirement that
$$ (x \cdot y, z) = ( y  \cdot z, x) $$
for all $x, y, x \in L$. This condition can be expressed equivalently by saying that $( u \cdot x , v) = (u, v \cdot x)$ for all $u, v, x \in L$ or in other words that multiplication from the right is self-adjoint with respect to the alternating form.\\ \\
Symplectic alternating algebras originate from a study of powerful 2-Engel groups~\cite{moravec},~\cite{gt-2008} and there is a 1-1 correspondence between a certain rich class of powerful 2-Engel 3-groups of exponent 27 and SAAs over the field $\mbox{GF}(3)$. We refer to~\cite{nsaa1},~\cite{sor-thesis} and~\cite{saa} for more detailed discussion of this as well as some general background to SAAs.\\ \\
In this paper we finish our classification of nilpotent symplectic alternating algebras of dimension 10 that we started in~\cite{nsaa2}.  Our approach like there relies on a general theory that we developed in~\cite{nsaa1}. Before giving some account of this we introduce some notation. Firstly we can always pick a basis $x_1, y_1, \ldots, x_n, y_n$ with the property that $(x_i, x_j)= (y_i, y_j)=0$ and $(x_i, y_j) = \delta_{ij}$ for $1 \leq i \leq j \leq k \leq n$. We refer to a basis of this type as \emph{standard basis}. It turns out that for any nilpotent symplectic alternating algebra one can always choose a suitable standard basis such that the chain of subspaces
\[0 = I_0 < I_1 < \cdots < I_n < I_{n-1}^\perp < \cdots < I_0^\perp = L, \]
with $I_k=\mathbb{F} x_n + \cdots + \mathbb{F} x_{n-k+1}$, is a central chain of ideals. One can furthermore see from this that $x_iy_j=0$ if $j \leq i$ and that $I_{n-1}^\perp$ is abelian. It follows that a number of the triple values $(uv, w)$ are trivial. Listing only the values that are possibly non-zero it suffices to consider
$${\mathcal P}: \quad (x_iy_j,y_k)=\alpha_{i j k},\quad (y_iy_j,y_k)=\beta_{i j k}, \quad 1 \leq i < j < k \leq n $$
for some $\alpha_{i j k}, \beta_{i j k} \in \mathbb{F}$. Such a presentation is called 
a \textit{nilpotent presentation}. 
Conversely any such presentations describes a nilpotent SAA. The algebras that are of maximal class turn out to have a rigid ideal structure. In particular when $2n \geq 10$ we can choose our chain of ideals above such that they are all characteristic and it turns out that $I_0, I_2, I_3, \ldots, I_{n-1}, I_{n-1}^\perp, I_{n-2}^\perp, \ldots, I_0^\perp$ are unique and equal to both the terms of the lower and upper central series ( see~\cite{nsaa1} Theorems 3.1 and 3.2 ). The algebras of maximal class can be identified easily from their nilpotent presentations. In fact, if ${\mathcal P}$ is any nilpotent presentation of $L$ with respect to a standard basis $\{ x_1, y_1, \ldots, x_n, y_n \}$, and $2n \geq 8$, we have that L is of maximal class if and only if $x_i y_{i+1} \neq 0$ for all $i = 2, \ldots, n-2$, and $x_1y_2, y_1y_2$ are linearly independent ( see~\cite{nsaa1} Theorem 3.4).
In~\cite{nsaa2} we found all the nilpotent algebras with a non-isotropic centre and we gave a full classification of the nilpotent SAAs with an isotropic centre of dimension 5 and 3. Here we deal with the two remaining cases, namely when the center has dimension 2 or 4. In the latter our approach leads to an interesting geometrical setup.
\section{Algebras with an isotropic centre of dimension 4}
The situation here is far more complicated and we will need to consider several
subcases. Let $L$ be a nilpotent SAA of dimension $10$ with an isotropic centre of
dimension $4$. We can pick our standard basis such that 
   $$Z(L)={\mathbb F}x_{5}+{\mathbb F}x_{4}+{\mathbb F}x_{3}+{\mathbb F}x_{2}.$$
Then 
     $$L^{2}=Z(L)^{\perp}=Z(L)+{\mathbb F}x_{1}+{\mathbb F}y_{1}$$
and by [2, Lemma 2.6] we know that $L^{3}=L^{2}L\leq Z(L)$. As $L^{2}\not\leq
Z(L)$ we have $L^{3}\not =\{0\}$ and by [2, Propositions 2.9 and 2.10] we then 
have that $3\leq \mbox{dim\,}L^{3}\leq 4$. We will consider the cases $L^{3}<Z(L)$
and $L^{3}=Z(L)$ separately. 
\subsection{The algebras where $L^{3}<Z(L)$}
In this case we can choose our basis such that 
$$\begin{array}{ll}
\begin{tabular}{c|c|c|c} 
\cline{2-3} 
   \mbox{} & $x_{5}$ & $y_{5}$ & \mbox{} \\
   \mbox{} & $x_{4}$ & $x_{4}$ & \mbox{} \\
   $L^{3}$ & $x_{3}$ & $y_{3}$ & \mbox{} \\
   \cline{2-3}  
  $Z(L)$ & $x_{2}$ & $y_{2}$ & $Z_{2}(L)$ \\
   \cline{2-3}
  \mbox{} & $x_{1}$ & $y_{1}$ & $L^{2}$ \\
\cline{2-3}
\end{tabular} &  
%
\begin{array}{l}
 L^{3}={\mathbb F}x_{5}+{\mathbb F}x_{4}+{\mathbb F}x_{3} \\ 
 Z(L)=L^{3}+{\mathbb F}x_{2} \\
 L^{2}=Z(L)^{\perp}=Z(L)+{\mathbb F}x_{1}+{\mathbb F}y_{1} \\
 Z_{2}(L)=(L^{3})^{\perp}=L^{2}+{\mathbb F}y_{2}
\end{array}
\end{array}$$
Notice that, as $Z_{2}(L)\cdot L^{2}=\{0\}$, we have that $Z_{2}(L)$
is abelian and thus in particular $x_{1}y_{2}=y_{1}y_{2}=0$. As $Z_{2}(L)$
is an abelian ideal we also have that $Z_{2}(L)L$ is orthogonal to $Z_{2}(L)$
and thus $Z_{2}(L)\cdot L\leq Z_{2}(L)^{\perp}=L^{3}$. It follows that 
   $${\mathbb F}x_{2}+{\mathbb F}x_{1}+{\mathbb F}y_{1}+L^{3}=
    L^{2}={\mathbb F}y_{3}y_{4}+{\mathbb F}y_{3}y_{5}+{\mathbb F}y_{4}y_{5}
     +L^{3}.$$
Suppose 
     $$x_{2}+L^{3}=\alpha y_{3}y_{4}+\beta y_{3}y_{5}+\gamma y_{4}y_{5}+L^{3}.$$
Now at least one of $\alpha,\beta,\gamma$ is nonzero and by the symmetry in
$y_{3},y_{4},y_{5}$ we can assume that $\alpha\not =0$. Thus 
     $$x_{2}+L^{3}=(y_{3}-\frac{\gamma}{\alpha} y_{5})(\alpha y_{4}+\beta y_{5})+L^{3}.$$
By replacing $x_{4},x_{5},y_{3},y_{4}$ by 
$\tilde{x}_{4}=\frac{1}{\alpha}x_{4}, 
\tilde{x}_{5}=x_{5}-\frac{\beta}{\alpha}x_{4}+\frac{\gamma}{\alpha}x_{3}, \tilde{y}_{3}=y_{3}-
\frac{\gamma}{\alpha}y_{5}, \tilde{y}_{4}=\alpha y_{4}+\beta y_{5}$, we can then assume that $x_{2}+L^{3}=y_{3}y_{4}+L^{3}$. In particular $(y_{2}y_{3},y_{4})=1$.
Suppose that $(y_{3}y_{4},y_{5})=\tau $. Replacing $x_{2},y_{5}$ by
$\tilde{x}_{2}=x_{2}+\tau x_{5}$ and $\tilde{y}_{5}=y_{5}-\tau y_{2}$ we can furthermore assume that $(y_{3}y_{4},y_{5})=0$. If we let $L={\mathbb F}y_{3}+
{\mathbb F}y_{4}+{\mathbb F}y_{5}$, it follows that we now have
\begin{equation}
    L^{2}={\mathbb F}x_{2}+{\mathbb F}x_{1}+{\mathbb F}y_{1}\mbox{\ and\ }
      y_{3}y_{4}=x_{2}.
\end{equation}
By (1) we know that 
     $$x_{1}=\alpha y_{3}y_{4}+\beta y_{3}y_{5}+\gamma y_{4}y_{5}=
           \alpha x_{2}+\beta y_{3}y_{5}+\gamma y_{4}y_{5},$$
where without loss of generality we can assume that $\gamma\not =0$. Then
     $$x_{1}=(y_{4}+\frac{\beta}{\gamma}y_{3})(-\alpha y_{3}+\gamma y_{5}).$$
Replacing $x_{3},x_{5},y_{4},y_{5}$ by $\tilde{x}_{3}=x_{3}-\frac{\beta}{\gamma}
x_{4}+\frac{\alpha}{\gamma}x_{5}$, $\tilde{x}_{5}=\frac{1}{\gamma}x_{5}$, 
$\tilde{y}_{4}=y_{4}+\frac{\beta}{\gamma}y_{3}$ and $\tilde{y}_{5}=-\alpha y_{3}+\gamma
y_{5}$, we obtain 
\begin{equation}
       y_{4}y_{5}=x_{1}.
\end{equation}
Notice that (1) is not affected by these changes. Finally we know from
(1) and (2) that 
     $$y_{1}=-\alpha x_{2}-\beta x_{1}+\gamma y_{5}y_{3}$$
for some $0 \not =\gamma\in {\mathbb F}$. Then 
    $$y_{1}=(y_{5}+\frac{\alpha}{\gamma} y_{4})(\gamma y_{3}+\beta y_{4}).$$
Now replace $x_{3},x_{4},y_{3},y_{5}$
by $\tilde{x}_{3}=\frac{1}{\gamma}x_{3}, \tilde{x}_{4}=x_{4}-
\frac{\alpha}{\gamma}x_{5}-\frac{\beta}{\gamma}x_{3}$, $\tilde{y}_{3}=\gamma y_{3}
+\beta y_{4}$ and $\tilde{y}_{5}=y_{5}+\frac{\alpha}{\gamma}y_{4}$. This gives us
\begin{equation}
       y_{5}y_{3}=y_{1}.
\end{equation}
This does not affect (2) but instead of (1) we get $y_{3}y_{4}=\gamma x_{2}$.
Now we make the final change by replacing $x_{2}$ and $y_{2}$
by $\gamma x_{2}$ and $\frac{1}{\gamma} y_{2}$ and we can assume that (1), (2)
and (3) hold. We had seen earlier that $Z_{2}(L)$ is abelian and thus all triple values involving two elements from $\{x_{5},x_{4},x_{3},x_{2},x_{1},y_{1},y_{2}\}$
is trivial. Thus all the nontrivial triple values involve two of $y_{3},y_{4},
y_{5}$ but from (1),(2) and (3) we know what these are. We have thus proved
\begin{prop} There is a unique nilpotent SAA of dimension $10$ that has an isotropic
center of dimension $4$ and where $L^{3}<Z(L)$. This algebra can be given
by the nilpotent presentation
  $${\mathcal P}_{10}^{(4,1)}:\ (y_{2}y_{3},y_{4})=1,\ (y_{1}y_{4},y_{5})=1,\ 
        (x_{1}y_{3},y_{5})=1.$$
\end{prop}
{\bf Remark}. Inspection shows that the algebra with that presentation has a centre
of dimension $4$ and the property that $L^{3}<Z(L)$.
\subsection{The algebras where $L^{3}=Z(L)$}
We will see that this case is quite intricate and we will need to 
consider some subcases.

$$\begin{array}{ll}
\begin{tabular}{c|c|c|c} 
\cline{2-3} 
   \mbox{} & $x_{5}$ & $y_{5}$ & \mbox{} \\
   \mbox{} & $x_{4}$ & $x_{4}$ & \mbox{} \\
   \mbox{} & $x_{3}$ & $y_{3}$ & \mbox{} \\  
  $Z(L)=L^{3}$ & $x_{2}$ & $y_{2}$ & \mbox{} \\
   \cline{2-3}
  \mbox{} & $x_{1}$ & $y_{1}$ & $L^{2}=Z_{2}(L)$ \\
\cline{2-3}
\end{tabular} &  
$$\begin{array}{l}
 Z(L)=L^{3}={\mathbb F}x_{5}+{\mathbb F}x_{4}+{\mathbb F}x_{3}+Fx_{2} \\ 
 Z_{2}(L)=L^{2}=Z(L)^{\perp}=Z(L)+{\mathbb F}x_{1}+{\mathbb F}y_{1} 
\end{array}$$
\end{array}$$
In order to clarify the structure further, we will associate to any such
algebra a family of new alternating forms that are defined as follows. For
each $\bar{z}=z+Z(L)\in L^{2}/Z(L)$, we obtain the alternating form
    $$\phi_{\bar{z}}:\,L/L^{2}\times L/L^{2}\rightarrow {\mathbb F}$$
given by 
    $$\phi_{\bar{z}}(\bar{u},\bar{v})=(zu,v)$$
where $\bar{u}=u+L^{2}$ and $\bar{v}=v+L^{2}$. Notice that this is a well defined function as $L^{2}$ is abelian. \\ \\
{\bf Remarks}. (1) If $0\not = \bar{z}=z+Z(L)\in L^{2}/Z(L)$, then $\phi_{\bar{z}}\not =0$. Otherwise we would have $(zu,v)=0$ for all $u,v\in L$ that would give
the contradiction that $z\in Z(L)$ and thus $\bar{z}=0$. \\ \\
(2) There is no non-zero element in $V=L/L^{2}$ that is common to the 
isotropic part of $V=L/L^{2}$ with respect to all the alternating
forms $\phi_{\bar{z}}$ with $\bar{z}\in L^{2}/Z(L)$.  
Otherwise there would be some $0\not =t\in {\mathbb F}y_{5}+
{\mathbb F}y_{4}+{\mathbb F}y_{3}+{\mathbb F}y_{2}$ such that $(zt,u)=0$ for
all $z\in L^{2}$ and all $u\in L$. But then $ut\in (L^{2})^{\perp}=Z(L)$ for
all $u\in L$ that gives the contradiction that $t\in Z_{2}(L)=L^{2}$. \\ \\
We divide the algebras into three 
categories. \\ \\
{\bf A}. The algebras where there exists a basis $\bar{z},\bar{t}$ for 
$L^{2}/Z(L)$ such that the alternating forms $\phi_{\bar{z}},\phi_{\bar{t}}$ are
both degenerate. \\ \\
{\bf B}. The algebras where there exists $0\not =\bar{z}\in L^{2}/Z(L)$ such that
$\phi_{\bar{z}}$ is degenerate but $\phi_{\bar{t}}$ is non-degenerate for
all $\bar{t}\in L^{2}/Z(L)$ that are not in ${\mathbb F}\bar{z}$. \\ \\
{\bf C}. The algebras where $\phi_{\bar{z}}$ is non-degenerate for
all $0\not =\bar{z}\in L^{2}/Z(L)$. 

\subsubsection{Algebras of type $A$}
Pick $x_{1}, y_{1}\in L^{2}\setminus Z(L)$ such that $\phi_{\bar{x}_{1}}$ and 
$\phi_{\bar{y}_{1}}$ are degenerate and such that $(x_{1},y_{1})=1$. 
By the remarks above we thus know that the isotropic part of $L/L^{2}$ with respect
to both the alternating forms $\phi_{\bar{x}_{1}}$ and $\phi_{\bar{y}_{1}}$ is
of dimension $2$ and the intersection of the two is trivial. Thus we can
pick a basis $\bar{y}_{5}=y_{5}+L^{2},\bar{y}_{4}=y_{4}+L^{2},
\bar{y}_{3}=y_{3}+L^{2},\bar{y}_{2}=y_{2}+L^{2}$ for $L/L^{2}$ such that 
$$\begin{array}{l}
      {\mathbb F}\bar{y}_{4}+{\mathbb F}\bar{y}_{5}\mbox{ is the isotropic
    part of }L/L^{2}\mbox{ with respect to }\phi_{\bar{x}_{1}}. \\
   {\mathbb F}\bar{y}_{3}+{\mathbb F}\bar{y}_{2}\mbox{ is the isotropic
    part of }L/L^{2}\mbox{ with respect to }\phi_{\bar{y}_{1}}. 
\end{array}$$
This shows us that we can pick our standard basis such that 
$$\begin{array}{ll}
   (x_{1}y_{2},y_{3})=1 & (y_{1}y_{2},y_{3})=0 \\
   (x_{1}y_{2},y_{4})=0 & (y_{1}y_{2},y_{4})=0 \\
   (x_{1}y_{2},y_{5})=0 & (y_{1}y_{2},y_{5})=0 \\
   (x_{1}y_{3},y_{4})=0 & (y_{1}y_{3},y_{4})=0 \\
   (x_{1}y_{3},y_{5})=0 & (y_{1}y_{3},y_{5})=0 \\
   (x_{1}y_{4},y_{5})=0 & (y_{1}y_{4},y_{5})=1.
\end{array}$$
To determine the structure fully we are only left with the triples $(y_{i}y_{j},y_{k})=r_{ijk}$ for $2\leq i<j<k\leq 5$. Let 
\begin{eqnarray*}
   \tilde{y}_{i} & = & y_{i}+\alpha_{i}x_{1}+\alpha_{i}y_{1}, \\ 
    \tilde{x}_{1} & = & x_{1}-(\alpha_{2}x_{2}+\alpha_{3}x_{3}+\alpha_{4}x_{4}
                    +\alpha_{5}x_{5}), \\
  \tilde{y}_{1} & = & y_{1}+\alpha_{2}x_{2}+\alpha_{3}x_{3}+\alpha_{4}x_{4}
                    +\alpha_{5}x_{5}.
\end{eqnarray*}
Inspection shows that we can choose $\alpha_{2},\ldots ,\alpha_{5}$ such that
$(\tilde{y}_{i}\tilde{y}_{j},\tilde{y}_{k})=0$ for all $2\leq i<j<k\leq 5$. In
fact this works for $\alpha_{2}=-r_{245}, \alpha_{3}=-r_{345}, \alpha_{4}=-r_{234}$
and $\alpha_{5}=-r_{235}$. We have thus proved the following result.
\begin{prop} There is a unique nilpotent SAA of dimension $10$ with an isotropic center 
of
dimension $4$ where $L^{3}=Z(L)$ and $L$ is of type A. This algebra
can be given by the presentation 
$${\mathcal P}_{10}^{(4,2)}:\ (x_{1}y_{2},y_{3})=1,\ (y_{1}y_{4},y_{5})=1.$$
\end{prop}
Notice that inspection shows that the algebra with this presentation indeed
has the properties stated in the proposition. 
\subsubsection{Algebras of type B}
Suppose that $\phi_{\bar{x}_{1}}$ is degenerate and that the isotropic part
of $L/L^{2}$ with respect to this alternating form is ${\mathbb F}\bar{y}_{4}
+{\mathbb F}\bar{y}_{5}$. We are now assuming that $\phi_{\bar{z}}$
is non-degenerate for all $\bar{z}\not\in {\mathbb F}\bar{x}_{1}$. Pick
$y_{1}\in L^{2}$ such that $(x_{1},y_{1})=1$. \\ \\
{\bf Remark}. We must have $\phi_{\bar{y}_{1}}(\bar{y}_{4},\bar{y}_{5})=0$. Otherwise we would get a basis $\bar{y}_{2},\bar{y}_{3},\bar{y}_{4},\bar{y}_{5}$ for 
$L/L^{2}$ such that 
   $$\phi_{\bar{y}_{1}}(\bar{y}_{4},\bar{y}_{5})=1,\ 
   \phi_{\bar{y}_{1}}(\bar{y}_{2},\bar{y}_{3})=1,\ 
   \phi_{\bar{x}_{1}}(\bar{y}_{2},\bar{y}_{3})=\alpha\not =0$$
and where $\phi_{\bar{y}_{1}}(\bar{y}_{i},\bar{y}_{j})=0$ and likewise 
$\phi_{\bar{x}_{1}}(\bar{y}_{i},\bar{y}_{j})=0$ for any pair $\bar{y}_{i},
\bar{y}_{j}$ such that $2\leq i<j\leq 5$ that is not included above. But then
inspection shows that $\phi_{\alpha \bar{y}_{1}-\bar{x}_{1}}$ is degenerate where
the corresponding isotropic part of $L/L^{2}$ is ${\mathbb F}\bar{y}_{2}+
{\mathbb F}\bar{y}_{3}$. But this contradicts our assumptions. \\ \\
We thus know that $\phi_{\bar{y}_{1}}(\bar{y}_{4},\bar{y}_{5})=0$. As 
$\phi_{\bar{y}_{1}}$ is non-degenerate we know that there exists some 
$\bar{y_{2}}\in L/L^{2}$ such that 
\begin{equation}
    \phi_{\bar{y}_{1}}(\bar{y}_{2},\bar{y}_{4})=1.
\end{equation} 
Replacing $y_{5}$ and $y_{3}$ by some suitable $y_{5}+\alpha y_{4}$ and
$y_{3}+\beta y_{4}+\gamma y_{2}$ we can furthermore assume that
\begin{equation}
 \phi_{\bar{y}_{1}}(\bar{y}_{2},\bar{y}_{5})=
 \phi_{\bar{y}_{1}}(\bar{y}_{3},\bar{y}_{4})=
 \phi_{\bar{y}_{1}}(\bar{y}_{2},\bar{y}_{3})=0. 
\end{equation}
As $\phi_{\bar{x}_{1}}$ is non-zero we must have $\phi_{\bar{x}_{1}}(\bar{y}_{2},
\bar{y}_{3})\not =0$ and by replacing $\bar{y_{3}}$ by a multiple of itself
we can assume that
\begin{equation}
    \phi_{\bar{x}_{1}}(\bar{y}_{2},\bar{y}_{3})=1.
\end{equation}
Notice that this does not affect (4) and (5). 
As $\phi_{\bar{y}_{1}}$ is non-degenerate we cannot have that $\bar{y}_{3}$ is
isotropic to all vectors in $L/L^{2}$ with respect to this alternating form.
Thus by (5) we must have $\phi_{\bar{y}_{1}}(\bar{y}_{3},\bar{y}_{5})\not =0$
and by replacing $\bar{y}_{5}$ by a multiple of itself we can assume that
\begin{equation}
  \phi_{\bar{y}_{1}}(\bar{y}_{3},\bar{y}_{5})=1.
\end{equation}
Again equations (4),(5) and (6) are not affected. We thus see that we can
choose a standard basis such that 
$$\begin{array}{ll}
   (x_{1}y_{2},y_{3})=1 & (y_{1}y_{2},y_{3})=0 \\
   (x_{1}y_{2},y_{4})=0 & (y_{1}y_{2},y_{4})=1 \\
   (x_{1}y_{2},y_{5})=0 & (y_{1}y_{2},y_{5})=0 \\
   (x_{1}y_{3},y_{4})=0 & (y_{1}y_{3},y_{4})=0 \\
   (x_{1}y_{3},y_{5})=0 & (y_{1}y_{3},y_{5})=1 \\
   (x_{1}y_{4},y_{5})=0 & (y_{1}y_{4},y_{5})=0.
\end{array}$$
As in case A we are now only left with the triples $(y_{i}y_{j},y_{k})=r_{ijk}$
for all $2\leq i<j<k\leq 5$. As in that case we let 
\begin{eqnarray*}
   \tilde{y}_{i} & = & y_{i}+\alpha_{i}x_{1}+\alpha_{i}y_{1}, \\ 
    \tilde{x}_{1} & = & x_{1}-(\alpha_{2}x_{2}+\alpha_{3}x_{3}+\alpha_{4}x_{4}
                    +\alpha_{5}x_{5}), \\
  \tilde{y}_{1} & = & y_{1}+\alpha_{2}x_{2}+\alpha_{3}x_{3}+\alpha_{4}x_{4}
                    +\alpha_{5}x_{5}.
\end{eqnarray*}
Inspection shows that we can choose $\alpha_{2},\alpha_{3},\alpha_{4},\alpha_{5}$ such that $(\tilde{y}_{i}\tilde{y}_{j},\tilde{y}_{k})=0$ for $2\leq i<j<k\leq 5$. We thus get the following result.
\begin{prop} There is a unique nilpotent SAA of dimension $10$ with isotropic center
of dimension $4$ where $L^{3}=Z(L)$ and $L$ is of type $B$. This algebra
can be given by the presentation 
   $${\mathcal P}_{10}^{(4,3)}:\,(x_{1}y_{2},y_{3})=1,\ (y_{1}y_{2},y_{4})=1,\ 
        (y_{1}y_{3},y_{5})=1.$$
\end{prop}
{\bf Proof}. We have already shown that this algebra is the only candidate. Inspection shows that conversely this algebra has isotropic centre of dimension
$4$ and $L^{3}=Z(L)$. It remains to see that the algebra is of type B. Thus
let $r\in {\mathbb F}$. We want to show that $\phi_{r\bar{x}_{1}+\bar{y}_{1}}$
is non-degenerate. Let $t=\alpha y_{2}+\beta y_{3}+\gamma y_{4}+
\delta y_{5}$ such that $\phi_{r\bar{x}_{1}+\bar{y}_{1}}(\bar{t},\bar{u})=0$
for all $\bar{u}\in L/L^{2}$ where $\bar{t}=t+L^{2}$. Then
$$\begin{array}{l}
   0=\phi_{r\bar{x}_{1}+\bar{y}_{1}}(\bar{t},\bar{y}_{5})=\beta \\
  0=\phi_{r\bar{x}_{1}+\bar{y}_{1}}(\bar{t},\bar{y}_{4})=\alpha \\
  0=\phi_{r\bar{x}_{1}+\bar{y}_{1}}(\bar{t},\bar{y}_{3})=r\alpha-\delta=-\delta \\
  0=\phi_{r\bar{x}_{1}+\bar{y}_{1}}(\bar{t},\bar{y}_{2})=-r\beta-\gamma=-\gamma.
\end{array}$$
Thus $\bar{t}=0$. $\Box$
\subsubsection{Algebras of type $C$}
Here we are assuming that $\phi_{z}$ is non-degenerate for all $0\not =
z\in L^{2}/Z(L)$.
Let $L$ be any nilpotent SAA of type $C$. Notice that $L^{2}/Z(L)=L^{2}/(L^{2})^{\perp}$
naturally becomes a $2$-dimensional symplectic vector space with inherited
alternating form from $L$. Thus $(u+Z(L),v+Z(L))=(u,v)$ for $u,v\in L^{2}$. 
We pick a basis $x,y$ for $L^{2}/Z(L)$ such that $(x,y)=1$ and then
choose some fixed elements $x_{1},y_{1}\in L^{2}$ such
$x=\bar{x}_{1}=x_{1}+Z(L)$ and $y=\bar{y}_{1}=y_{1}+Z(L)$. 
For any
vector $u\in L/L^{2}$ we will denote by $\langle u\rangle_{1}^{\perp}$ the
subspace of $L/L^{2}$ consisting of all the vectors that are isotropic
to $u$ with respect to $\phi_{\bar{x}_{1}}$. Likewise we will denote
by $\langle u\rangle_{2}^{\perp}$ the subspace of $L/L^{2}$ consisting of
all the vectors that are isotropic to $u$ with respect to 
$\phi_{\bar{y}_{1}}$. \\ \\
{\bf Definition}. We say that a subspace of $L/L^{2}$ is {\it totally isotropic}
if it is isotropic with respect to $\phi_{z}$ for all $z\in 
L^{2}/Z(L)$. 
\begin{lemm} For each $0\not =u\in L/L^{2}$
there exists a unique totally isotropic plane through $u$. 
\end{lemm} 
{\bf Proof}\ \ We know that $\langle u\rangle_{1}^{\perp}$ and $\langle u
\rangle_{2}^{\perp}$ are $3$-dimensional. Thus if they are not equal then 
    $$4=\mbox{dim\,}(\langle u\rangle_{1}^{\perp}+\langle u\rangle_{2}^{\perp})
      =\mbox{dim\,}\langle u\rangle_{1}^{\perp}+\mbox{dim\,}\langle
u\rangle_{2}^{\perp}-\mbox{dim\,}(\langle u\rangle_{1}^{\perp}\cap 
\langle u\rangle_{2}^{\perp}).$$
Therefore $\mbox{dim\,}(\langle u\rangle_{1}^{\perp}\cap 
\langle u\rangle_{2}^{\perp})=6-4=2$. 
Thus the collection of all the elements in $L/L^{2}$ that are isotropic to
$u$ with respect to $\phi_{z}$ for all $z\in L^{2}/Z(L)$, namely
$\langle u\rangle_{1}^{\perp}\cap \langle u\rangle_{2}^{\perp}$, is a 
plane. \\ \\
It remains to see that $\langle u\rangle_{1}^{\perp}\not =
\langle u\rangle_{2}^{\perp}$. We argue by contradiction and pick a basis
$u,v,w$ for this common subspace and add a vector $t$ to get a basis
for $L/L^{2}$. By replacing $x_{1}$ by some suitable $x_{1}+
\alpha y_{1}$, we can assume that $\phi_{\bar{x}_{1}}(u,t)=0$. But
then $u$ is isotropic to all elements of $L/L^{2}$ with respect to
$\phi_{\bar{x}_{1}}$ that contradicts the assumption that $\phi_{\bar{x}_{1}}$
is non-degenerate. $\Box$ \\ \\
The alternating forms $\phi_{z}$ with $z\in L^{2}/Z(L)$ will help
us understanding the structure of algebras of type C. 
We will next come up with a special type of presentations for algebras of type C based on the geometry arising from the family of the auxiliary alternating 
forms. \\ \\
Let $L$ be an algebra of type C. As a starting 
point we pick two distinct totally isotropic planes $P_{1},P_{2}\leq L/L^{2}$
and we pick some non-zero vector $\bar{y}_{2}$ on $P_{1}$. By Lemma 2.4, we have that
$P_{1}\cap P_{2}=\{0\}$ and thus $L/L^{2}=P_{1}\oplus P_{2}$. Now $\langle \bar{y}_{2}\rangle_{1}^{\perp}$ is $3$-dimensional and contains $P_{1}$. Thus $\langle 
\bar{y}_{2}\rangle_{1}^{\perp}\cap P_{2}$ is $1$-dimensional and not contained
in $\langle \bar{y}_{2}\rangle_{2}^{\perp}$ by Lemma 2.4. Thus there is unique vector
$\bar{y}_{5}\in P_{2}$ where $\phi_{\bar{x}_{1}}(\bar{y}_{2},\bar{y}_{5})=0$ and 
$\phi_{\bar{y}_{1}}(\bar{y}_{2},\bar{y}_{5})=1$. 
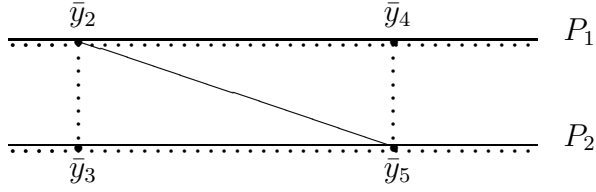
\mbox{}\\
\begin{picture}(200,90)(-100,30)
\put(-25,100){\line(3,-1){120}}

\put(96,99){$\bulletlike$}
\put(96,59){$\bulletlike$}
\put(92,107){$\bar{y}_{4}$}
\put(92,48){$\bar{y}_{5}$}
\put(-50,100){\line(1,0){200}}
\put(160,100){$P_{1}$}
\put(160,60){$P_{2}$}
\put(-25,95){$\cdot$}
\put(-25,90){$\cdot$}
\put(-25,85){$\cdot$}
\put(-25,80){$\cdot$}
\put(-25,75){$\cdot$}
\put(-25,70){$\cdot$}
\put(-25,65){$\cdot$}
\put(-25,60){$\cdot$}
\put(94,95){$\cdot$}
\put(94,90){$\cdot$}
\put(94,85){$\cdot$}
\put(94,80){$\cdot$}
\put(94,75){$\cdot$}
\put(94,70){$\cdot$}
\put(94,65){$\cdot$}
\put(94,60){$\cdot$}
%

\put(-23.5,99){$\bulletlike$}
\put(-27,107){$\bar{y}_{2}$}
\put(-27,48){$\bar{y}_{3}$}
\put(-23.5,59){$\bulletlike$}
\put(-50,95){$\cdot$}
\put(-45,95){$\cdot$}
\put(-40,95){$\cdot$}
\put(-35,95){$\cdot$}
\put(-50,95){$\cdot$}
\put(-45,95){$\cdot$}
\put(-40,95){$\cdot$}
\put(-35,95){$\cdot$}
\put(-30,95){$\cdot$}
\put(-25,95){$\cdot$}
\put(-20,95){$\cdot$}
\put(-15,95){$\cdot$}
\put(-10,95){$\cdot$}
\put(-5,95){$\cdot$}
\put(0,95){$\cdot$}
\put(5,95){$\cdot$}
\put(10,95){$\cdot$}
\put(15,95){$\cdot$}
\put(20,95){$\cdot$}
\put(25,95){$\cdot$}
\put(30,95){$\cdot$}
\put(35,95){$\cdot$}
\put(40,95){$\cdot$}
\put(45,95){$\cdot$}
\put(50,95){$\cdot$}
\put(55,95){$\cdot$}
\put(60,95){$\cdot$}
\put(65,95){$\cdot$}
\put(70,95){$\cdot$}
\put(75,95){$\cdot$}
\put(80,95){$\cdot$}
\put(85,95){$\cdot$}
\put(90,95){$\cdot$}
\put(95,95){$\cdot$}
\put(100,95){$\cdot$}
\put(105,95){$\cdot$}
\put(110,95){$\cdot$}
\put(115,95){$\cdot$}
\put(120,95){$\cdot$}
\put(125,95){$\cdot$}
\put(130,95){$\cdot$}
\put(135,95){$\cdot$}
\put(140,95){$\cdot$}
\put(145,95){$\cdot$}
\put(-50,60){\line(1,0){200}}
\put(-50,55){$\cdot$}
\put(-45,55){$\cdot$}
\put(-40,55){$\cdot$}
\put(-35,55){$\cdot$}
\put(-50,55){$\cdot$}
\put(-45,55){$\cdot$}
\put(-40,55){$\cdot$}
\put(-35,55){$\cdot$}
\put(-30,55){$\cdot$}
\put(-25,55){$\cdot$}
\put(-20,55){$\cdot$}
\put(-15,55){$\cdot$}
\put(-10,55){$\cdot$}
\put(-5,55){$\cdot$}
\put(0,55){$\cdot$}
\put(5,55){$\cdot$}
\put(10,55){$\cdot$}
\put(15,55){$\cdot$}
\put(20,55){$\cdot$}
\put(25,55){$\cdot$}
\put(30,55){$\cdot$}
\put(35,55){$\cdot$}
\put(40,55){$\cdot$}
\put(45,55){$\cdot$}
\put(50,55){$\cdot$}
\put(55,55){$\cdot$}
\put(60,55){$\cdot$}
\put(65,55){$\cdot$}
\put(70,55){$\cdot$}
\put(75,55){$\cdot$}
\put(80,55){$\cdot$}
\put(85,55){$\cdot$}
\put(90,55){$\cdot$}
\put(95,55){$\cdot$}
\put(100,55){$\cdot$}
\put(105,55){$\cdot$}
\put(110,55){$\cdot$}
\put(115,55){$\cdot$}
\put(120,55){$\cdot$}
\put(125,55){$\cdot$}
\put(130,55){$\cdot$}
\put(135,55){$\cdot$}
\put(140,55){$\cdot$}
\put(145,55){$\cdot$}
\end{picture}
\normalsize
\mbox{}\\
Similarly we have a unique element $\bar{y}_{3}\in P_{2}$ such that 
$\phi_{\bar{y}_{1}}(\bar{y}_{2},\bar{y}_{3})=0$ and  $\phi_{\bar{x}_{1}}(\bar{y}_{2},\bar{y}_{3})=1$. By Lemma 2.4 we have 
$\langle \bar{y}_{5}\rangle_{2}^{\perp}\not =\langle \bar{y}_{3}
\rangle_{2}^{\perp}$. Thus there exists a unique $\bar{y}_{4}\in P_{1}$ such that
$\phi_{\bar{y}_{1}}(\bar{y}_{4},\bar{y}_{5})=0$
and $\phi_{\bar{y}_{1}}(\bar{y}_{3},\bar{y}_{4})=1$. Notice also that $\phi_{\bar{x}_{1}}(\bar{y}_{4},
\bar{y}_{5})\not =0$ and that $\bar{y}_{2},\bar{y}_{3},\bar{y}_{4},\bar{y}_{5}$
form a basis for $L/L^{2}$. It follows from the discussion that, for some $\alpha,\beta \in {\mathbb F}$ with $\beta \not =0$, we have 
\begin{equation}\begin{array}{llll}
\phi_{\bar{x}_{1}}(\bar{y}_{2},\bar{y}_{3})=1, & \phi_{\bar{y}_{1}}(\bar{y}_{2},\bar{y}_{3})=0, & \phi_{\bar{x}_{1}}(\bar{y}_{2},\bar{y}_{4})=0, &
\phi_{\bar{y}_{1}}(\bar{y}_{2},\bar{y}_{4})=0, \\
\phi_{\bar{x}_{1}}(\bar{y}_{2},\bar{y}_{5})=0, & \phi_{\bar{y}_{1}}(\bar{y}_{2},\bar{y}_{5})=1, & \phi_{\bar{x}_{1}}(\bar{y}_{3},\bar{y}_{4})=\alpha, &
\phi_{\bar{y}_{1}}(\bar{y}_{3},\bar{y}_{4})=1, \\
\phi_{\bar{x}_{1}}(\bar{y}_{3},\bar{y}_{5})=0, & \phi_{\bar{y}_{1}}(\bar{y}_{3},\bar{y}_{5})=0, & \phi_{\bar{x}_{1}}(\bar{y}_{4},\bar{y}_{5})=\beta, &
\phi_{\bar{y}_{1}}(\bar{y}_{4},\bar{y}_{5})=0.
\end{array}
\end{equation}
The matrix for the alternating form $\phi_{r\bar{x}_{1}+s\bar{y}_{1}}$ with
respect to the ordered  basis $(\bar{y}_{2},\bar{y}_{4},\bar{y}_{3},
\bar{y}_{5})$ is then
$$r\left[\begin{array}{rrrr}
     0 & 0 & 1 & 0 \\
     0 & 0 & -\alpha & \beta \\
    -1 & \alpha & 0 & 0 \\
    0 & -\beta & 0 & 0
\end{array}\right]+
s\left[\begin{array}{rrrr}
     0 & 0 & 0 & 1 \\
     0 & 0 & -1 & 0 \\
    0 & 1 & 0 & 0 \\
    -1 & 0 & 0 & 0
\end{array}\right]
$$
with determinant $(\beta r^{2}+\alpha rs+s^{2})^{2}$. As we are dealing here
with algebras of type C this determinant must be non-zero for all 
$(r,s)\not =(0,0)$. Equivalently we must have that the polynomial 
      $$t^{2}+\alpha t+\beta$$
is irreducible in ${\mathbb F}[t]$. Using this and (8) we will now obtain
a full presentation for our algebra. As before we are only left with the triples $(y_i y_j, y_k) = r_{ijk}$ for $2 \leq i < j < k \leq 5$. We will see that we can choose a standard 
basis such that $x_{1}+Z(L)=x$, $y_{1}+Z(L)=y$ and
$y_{i}+L^{2}=\bar{y}_{i}$ for $2\leq i\leq 5$. It turns out that we do not
have to alter our basis elements $x_{5},\ldots ,x_{2}$ of the centre. We do
this with a change of basis of the form 
\begin{eqnarray*}
  \tilde{x}_{1} & = & x_{1}-(\alpha_{2}x_{2}+\alpha_{3}x_{3}+\alpha_{4}x_{4}+\alpha_{5}x_{5}) \\
\tilde{y}_{1} & = & y_{1}+ (\alpha_{2}x_{2}+\alpha_{3}x_{3}+\alpha_{4}x_{4}+\alpha_{5}x_{5}) \\
 \tilde{y}_{i} & = & y_{i}+\alpha_{i}x_{1}+\alpha_{i}y_{1}.
\end{eqnarray*}
Inspection shows that the equations 
$(\tilde{y}_{i}\tilde{y}_{j},\tilde{y}_{k})=0$, $2\leq i<j<k\leq 5$ are 
equivalent to 
\begin{eqnarray*}
\left[\begin{array}{cc}
    -1 & 1 \\
    \beta & \alpha +1 
\end{array}\right]
\left[\begin{array}{r}
   \alpha_{3} \\
   \alpha_{5} 
\end{array}\right] &
= & \left[
\begin{array}{r}
   -r_{235} \\
   -r_{345}
\end{array}\right], \\
\left[\begin{array}{cc}
    \alpha +1  & 1 \\
    \beta & -1 
\end{array}\right]
\left[\begin{array}{r}
   \alpha_{2} \\
   \alpha_{4} 
\end{array}\right] &
= & \left[
\begin{array}{r}
   -r_{234} \\
   -r_{245}
\end{array}\right].
\end{eqnarray*}
Notice that we cannot have that $\alpha +\beta +1=0$ since otherwise $1$ is
a root of $t^{2}+\alpha t+\beta$ that is absurd as the polynomial
is irreducible. We thus have solution $(\alpha_{2},\alpha_{3},\alpha_{4},
\alpha_{5})$ to the equation system and we arrive at a standard basis
that gives us the following presentation
$${\mathcal P}_{10}^{(4,4)}(\alpha,\beta):\ (x_{1}y_{2},y_{3})=1,\ 
 (x_{1}y_{3},y_{4})=\alpha,\ (x_{1}y_{4},y_{5})=\beta,\ 
 (y_{1}y_{2},y_{5})=1,\ (y_{1}y_{3},y_{4})=1,$$
where the polynomial $t^{2}+\alpha t +\beta$ is irreducible. Conversely, inspection shows that any algebra with such presentation where $t^{2}+\alpha t+\beta$
is irreducible, gives us an algebra of type C. \\ \\
We next turn to the isomorphism problem, that is we want to understand when
two pairs $(\alpha,\beta)$ and $(\tilde{\alpha},\tilde{\beta})$ describe
the same algebra. As a first step we first prove the following lemma.
\begin{lemm} Let $x,y$ be elements in $L^{2}/Z(L)$
such that $(x,y)=1$. We have that the values of $\alpha$ and $\beta$ remain the same for all presentations of the form 
${\mathcal P}_{10}^{(4,4)}(\alpha,\beta)$ where, for the given standard basis, $x_{1}+Z(L)
=x$ and $y_{1}+Z(L)=y$. 
\end{lemm}
{\bf Proof}\ \ Our method for producing $\alpha$ and $\beta$ was based on choosing some distinct totally isotropic planes $P_{1},P_{2}$ and some non-zero vector 
$\bar{y}_{2}$ on $P_{1}$. From this we came up with a procedure that provided
us with unique vectors $\bar{y}_{2},\bar{y}_{4}\in P_{1}$ and 
$\bar{y}_{3},\bar{y}_{5}\in P_{2}$ from which the values $\alpha$ and $\beta$ can be
calculated as 
   $$\alpha = \phi_{x}(\bar{y}_{3},\bar{y}_{4}),\ \beta=
     \phi_{y}(\bar{y}_{4},\bar{y}_{5}).$$
We want to show that if $x_{1}+Z(L),y_{1}+Z(L)$ are kept fixed the procedure 
will always produce the same value for $\alpha$ and $\beta$. As a starting
point we will see that the values do not depend on which vector $\bar{y}_{2}$
from $P_{1}$ we choose. We do this in two steps. First notice that if we choose
$\tilde{y}_{2}=a\bar{y}_{2}$ for some $0\not =a\in {\mathbb F}$, then 
the procedure gives us the new vectors $\tilde{y}_{3}=\frac{1}{a}\bar{y}_{3}$,
$\tilde{y}_{5} = \frac{1}{a}\bar{y}_{5}$ and $\tilde{y}_{4}=ay_{4}$ and this gives us
the values 
$$\begin{array}{l}
 \tilde{\alpha}=\phi_{\bar{x}_{1}}(\tilde{y}_{3},\tilde{y}_{4})=
   \phi_{\bar{x}_{1}}(\frac{1}{a}\bar{y}_{3},a\bar{y}_{4})=\alpha, \\
\tilde{\beta}=\phi_{\bar{x}_{1}}(\tilde{y}_{4},\tilde{y}_{5})=
        \phi_{\bar{x}_{1}}(a\bar{y}_{4},\frac{1}{a}\bar{y}_{5})=\beta.
\end{array}$$
It remains to consider the change $\tilde{y}_{2}=\bar{y}_{4}+a\bar{y}_{2}$. Following the mechanical procedure above produces the elements 
\begin{eqnarray*}
   \tilde{y}_{5} & = & \frac{-\beta}{a^{2}-a\alpha+\beta}\,\bar{y}_{3}+
        \frac{a-\alpha}{a^{2}-a\alpha+\beta}\,\bar{y}_{5}, \\
   \tilde{y}_{3} & = & \frac{a}{a^{2}-a\alpha+\beta}\,\bar{y}_{3}+
        \frac{1}{a^{2}-a\alpha+\beta}\,\bar{y}_{5}, \\
    \tilde{y}_{4} & = &-\beta\,\bar{y}_{2}+
        (a-\alpha)\,\bar{y}_{4}.
\end{eqnarray*}
Inspection shows that again we have $\phi_{\bar{x}_{1}}(\tilde{y}_{3},
\tilde{y}_{4})=\alpha$ and $\phi_{\bar{x}_{1}}(\tilde{y}_{4},
\tilde{y}_{5})=\beta$. \\ \\
We have thus seen that for a given pair $P_{1},P_{2}$ of distinct totally
isotropic planes we get unique values $\alpha(P_{1},P_{2})$ and
$\beta(P_{1},P_{2})$ not depending on which vector $\bar{y}_{2}$ from
$P_{1}$ we choose for the procedure above. The next step is to see
that $\alpha(P_{2},P_{1})=\alpha(P_{1},P_{2})$ and $\beta(P_{2},P_{1})=
\beta(P_{1},P_{2})$. So suppose we have some standard basis with respect
to the pair $P_{1},P_{2}$ that gives us the presentation ${\mathcal P}_{10}^{(4,4)}(\alpha,
\beta)$. Now consider $\tilde{y}_{2}=\bar{y}_{5},\tilde{y}_{4}=\beta
\bar{y}_{3}\in P_{2}$ and $\tilde{y}_{3}=\frac{-1}{\beta}\bar{y}_{4},
\tilde{y}_{5}=-\bar{y}_{2}\in P_{1}$. Inspection shows that this is a standard
basis for the pair $P_{2},P_{1}$ and 
$$\begin{array}{l}
  \alpha(P_{2},P_{1})=\phi_{\bar{x}_{1}}(\tilde{y}_{3},\tilde{y}_{4})=
       \phi_{\bar{x}_{1}}(\frac{-1}{\beta}\bar{y}_{4},\beta \bar{y}_{3})=
\alpha(P_{1},P_{2}) \\
  \beta(P_{2},P_{1})=\phi_{\bar{x}_{1}}(\tilde{y}_{4},\tilde{y}_{5})=
       \phi_{\bar{x}_{1}}(\beta \bar{y}_{3},-\bar{y}_{2})=\beta(P_{1},P_{2}).
\end{array}$$
Now pick any totally isotropic plane $P_{3}$ that is distinct from $P_{1},
P_{2}$. The aim is to show that $\alpha(P_{3},P_{2})=\alpha(P_{1},P_{2})$
and $\beta(P_{3},P_{2})=\beta(P_{1},P_{2})$. Take any basis for $P_{3}$. This 
must be of the form $u_{1}+v_{1},u_{2}+v_{2}$ with $u_{1},u_{2}\in P_{1}$
and $v_{1},v_{2}\in P_{2}$. Notice first that $u_{1}+P_{2},u_{2}+P_{2}$ are linearly independent vectors in $P_{1}+P_{2}/P_{2}$. To see this, take $a,b\in {\mathbb F}$ such that 
    $$P_{2}=a u_{1}+b u_{2}+P_{2}=a (u_{1}+v_{1})+b (u_{2}+v_{2})
+P_{2}.$$
Then $a (u_{1}+v_{1})+b (u_{2}+v_{2})\in P_{2}\cap P_{3}=\{0\}$. As the
vectors $u_{1}+v_{1},u_{2}+v_{2}$ are linearly independent it follows that
$a=b=0$. In particular we can choose our basis for $P_{3}$ to be of the form
$\bar{y}_{2}+u, \bar{y}_{4}+v$ with $u,v\in P_{2}$. Inspection shows that for
$\tilde{y}_{2}=\bar{y}_{2}+u,\tilde{y}_{4}=\bar{y}_{4}+v\in P_{3}$ and
$\bar{y}_{3},\bar{y}_{5}\in P_{2}$ we have a standard basis with respect
to the pair $P_{3},P_{2}$. Furthermore the corresponding parameters are
$\phi_{\bar{x}_{1}}(\bar{y}_{3},\tilde{y}_{4})=\alpha$ and $\phi_{\bar{x}_{1}}(
\tilde{y}_{4},\bar{y}_{5})=\beta$.\\ \\
We have now all the input we need to finish the proof of the lemma. Take 
any four totally isotropic planes $P_{1},P_{2},P_{3},P_{4}$ in $L/L^{2}$ such
that $P_{1}\not =P_{2}$ and $P_{3}\not =P_{4}$. If these planes are not
all distinct then we get directly from the analysis above that $\alpha(P_{3},P_{4})=\alpha(P_{1},P_{2})$ and $\beta(P_{3},P_{4})=\beta(P_{1},P_{2})$. Now suppose
the planes are distinct. Then $\alpha(P_{3},P_{4})=\alpha(P_{1},P_{4})=
\alpha(P_{1},P_{2})$ and 
$\beta(P_{3},P_{4})=\beta(P_{1},P_{4})=\beta(P_{1},P_{2})$. This finishes
the proof of the lemma. $\Box$ \\ \\
It follows from the lemma that if we want to obtain a new presentation for
some given algebra $L$, then we must choose different vectors $x,y$ for $L^{2}/Z(L)$. We thus only need to consider 
a change of standard basis for $L$ of the form  $\tilde{x}_{1}=ax_{1}+by_{1},
\tilde{y}_{1}=cx_{1}+dy_{1}$ where $1=(\tilde{x}_{1},\tilde{y}_{1})=ad-bc$.
Suppose that we have a presentation ${\mathcal P}_{10}^{(4,4)}(\alpha,\beta)$ with respect
to some standard basis $x_{1},y_{1},\ldots ,x_{5},y_{5}$ and let $\tilde{x}_{1},
\tilde{y}_{1}$ be given as above. Going again through the standard procedure
with respect to $P_{1}={\mathbb F}\bar{y}_{2}+{\mathbb F}\bar{y}_{4}$,
$P_{2}={\mathbb F}\bar{y}_{3}+{\mathbb F}\bar{y}_{5}$ and $\bar{y}_{2}\in P_{1}$
gives us the new basis $\bar{y}_{2},\tilde{y}_{3},\tilde{y}_{5},\tilde{y}_{4}$
where
\begin{eqnarray*}
  \tilde{y}_{5} & = & -b\bar{y}_{3}+ a\bar{y}_{5} \\
  \tilde{y}_{3} & = & d\bar{y}_{3}-c\bar{y}_{5} \\
  \tilde{y}_{4} & = & 
\frac{-\alpha bc-\beta ac-bd}{d^{2}+cd \alpha +\beta c^{2}}\,\bar{y}_{2}+
\frac{1}{d^{2}+\alpha cd+\beta c^{2}}\,\bar{y}_{4}.
\end{eqnarray*}
From this we can calculate the new parameters $\tilde{\alpha}$ and
$\tilde{\beta}$ and we obtain the following proposition.
\begin{prop} Let $L$ be a nilpotent SAA of dimension $10$ with an isotropic centre of dimension
$4$ where $L^3=Z(L)$ and $L$ is of type C. Then $L$ has a presentation of the form
$${\mathcal P}_{10}^{(4,4)}(\alpha,\beta):\ (x_{1}y_{2},y_{3})=1,\ 
 (x_{1}y_{3},y_{4})=\alpha,\ (x_{1}y_{4},y_{5})=\beta,\ 
 (y_{1}y_{2},y_{5})=1,\ (y_{1}y_{3},y_{4})=1,$$
where the polynomial $t^{2}+\alpha t+\beta$ is irreducible in ${\mathbb F}[t]$.
Conversely any such presentation gives us an algebra of type C. Furthermore
two presentations ${\mathcal P}_{10}^{(4,4)}(\alpha,\beta)$ and 
${\mathcal P}_{10}^{(4,4)}(\tilde{\alpha},
\tilde{\beta})$ describe the same algebra if and only if 
\begin{eqnarray*}
\tilde{\alpha} & = & \frac{(ad+bc)\alpha+2ac\beta +2bd}{d^{2}+\alpha cd+
c^{2}\beta} \\
\tilde{\beta} & = & \frac{b^{2}+ab\alpha +a^{2}\beta}{d^{2}+cd \alpha +
c^{2}\beta}
\end{eqnarray*}
for some $a,b,c,d\in {\mathbb F}$ where $ad-bc=1$. 
\end{prop}
\subsection{Further analysis of algebras of  type $C$ and some examples}
In order to get a more transparent picture of the algebras of type C, it turns
out to be useful to consider the case when the characteristic is $2$ and the
case when the characteristic is not $2$ separately. 
\begin{lemm} Let $L$ be a nilpotent SAA of dimension $10$ with an isotropic center of dimension $4$ where $L^3=Z(L)$ and $L$ is of type C over a field ${\mathbb F}$ of 
characteristic that is not $2$. Then $L$ has a presentation of the form
${\mathcal P}_{10}^{(4,4)}(0,\beta)$ with respect to some standard basis, where the polynoimal $t^2+\beta$ is irreducible in $\mathbb{F}[t]$.
\end{lemm}
{\bf Proof}. By Proposition 2.6 we know that we can choose a standard basis
for $L$ so that $L$ has a presentation of the form ${\mathcal P}_{10}^{(4,4)}(\alpha,\beta)$
with respect to this basis. Now let $a=0$, $b=1$, $c=-1$ and $d=\alpha/2$. Then
$ad-bc=1$ and by Proposition 2.6 again we know that there is a presentation
for $L$ of the form ${\mathcal P}_{10}^{(4,4)}(\tilde{\alpha},\tilde{\beta})$ where
$\tilde{\alpha}=0$. $\Box$
\begin{prop} Let $L$ be a nilpotent SAA of dimension $10$ with an isotropic center of dimension $4$ where $L^3=Z(L)$ and $L$ is of type C over a field ${\mathbb F}$ of characteristic that is not $2$. Then $L$ has a presentation of the form 
$${\mathcal P}(\beta):\ (x_{1}y_{2},y_{3})=1,\ 
 (x_{1}y_{4},y_{5})=\beta,\ 
 (y_{1}y_{2},y_{5})=1,\ (y_{1}y_{3},y_{4})=1,$$
where $\beta \not\in -{\mathbb F}^{2}$. Conversely any such presentation gives 
us an algebra of type C. Furthermore two such presentations ${\mathcal P}(\beta)$ and ${\mathcal P}(\tilde{\beta})$ describe the same algebra if and only if
    $$\tilde{\beta}/\beta=(a^{2}+b^{2}\beta)^{2}$$
for some $(a,b)\in {\mathbb F}\times {\mathbb F}\setminus \{(0,0)\}$. 
\end{prop}
{\bf Proof}\ \ From Lemma 2.7 we know that such a presentation exists and the 
polynomial $t^{2}+\beta$ is irreducible if and only if 
$\beta\not\in -{\mathbb F}^{2}$. By Proposition 2.6 we then know that any such presentation gives us an algebra of type C and that
${\mathcal P}(0,\beta)$ and ${\mathcal P}(0,\tilde{\beta})$ describe the
same algebra if and only if there are $a,b,c,d\in {\mathbb F}$ such that
   $$0=ac\beta+bd,\ \ ad-bc=1$$
and 
     $$\tilde{\beta}=\frac{b^{2}+
   a^{2}\beta}{d^{2}+c^{2}\beta}.$$
Solving these together shows that these conditions are equivalent to
saying that 
   $$\tilde{\beta}=((\frac{b}{\beta})^{2}\beta + a^{2})^{2}\beta$$
for some $(a,b) \in {\mathbb F} \times {\mathbb F} \setminus \{ (0,0)\} $. As $\frac{b}{\beta}$ is arbitrary, the
second part of the proposition follows. $\Box$
\mbox{}\\ \\
{\bf Examples}. (1) If ${\mathbb F}={\mathbb C}$ then there are no algebras
of type C. This holds more generally for any field ${\mathbb F}$ whose characteristic is not $2$ and where all elements in ${\mathbb F}$ have a square root
in ${\mathbb F}$. \\ \\
(2) Suppose ${\mathbb F}={\mathbb R}$. Here $\beta \not \in -{\mathbb R}^{2}$
if and only if $\beta > 0$ in which case there exist $a \in {\mathbb R} \setminus  \{0\}$
such that $\beta\, =\, 1/a^{4}$. Hence, by Proposition 2.8,  ${\mathcal P}(\beta)$
describes the same algebra as ${\mathcal P}(1)$. This shows that there is only
one algebra of type C over ${\mathbb R}$ that can be given by the presentation
 $${\mathcal P}(1):\ (x_{1}y_{2},y_{3})=1,\ 
 (x_{1}y_{4},y_{5})=1,\ 
 (y_{1}y_{2},y_{5})=1,\ (y_{1}y_{3},y_{4})=1.$$
(3) Let ${\mathbb F}$ be a finite field of some odd characteristic $p$. Suppose
that $|{\mathbb F}|=p^{n}$. The nonzero elements form a cyclic group 
${\mathbb F}^{*}$ of order $p^{n}-1$ that is divisible by $2$. Thus there
are two cosets of $({\mathbb F}^{*})^{2}$ in ${\mathbb F}^{*}$ and 
   $${\mathbb F}^{*}=-({\mathbb F}^{*})^{2}\cup \beta({\mathbb F}^{*})^{2}$$
for some $\beta\in {\mathbb F}^{*}$. 
Suppose $\tilde{\beta}=\beta c^{2}$ is an arbitrary field element that
is not in $ - {\mathbb F}^{2}$. 
As the are $(|{\mathbb F}|+1)/2$ squares in ${\mathbb F}$ we have that the set $c - {\mathbb F}^{2}\ $ and $\beta {\mathbb F}^{2}$ intersect. Hence there exist
$a,b\in {\mathbb F}$ such that $c=a^{2}+b^{2}\beta$ and thus
$\tilde{\beta}=(a^{2}+b^{2}\beta)^{2}\beta$. Hence the situation is like
in (2) and we get only one algebra with presentation 
  $${\mathcal P}(\beta):\ (x_{1}y_{2},y_{3})=1,\ 
 (x_{1}y_{4},y_{5})=\beta,\ (y_{1}y_{2},y_{5})=1,\ (y_{1}y_{3},y_{4})=1,$$
where $\beta$ is any element not in $-{\mathbb F}^{2}$. \\ \\
{\bf Remark}. Let $\beta\in {\mathbb F}$ that is not in $-{\mathbb F}^{2}$
and consider a splitting field ${\mathbb F}[\gamma ]$ of the polynomial 
$t^{2}+\beta$ in ${\mathbb F}[t]$ where $\gamma^{2}=-\beta$. So $a^{2}+b^{2}\beta$ is
the norm $N(a+\gamma b)$ of $a+\gamma b$ that is a multiplicative function and thus 
     $$G(\beta)=\{(a^{2}+b^{2}\beta)^{2} :\, (a,b)\in {\mathbb F}\times
{\mathbb F}\setminus \{(0,0)\}\}$$
is a subgroup of $({\mathbb F}^{*})^{2}$. Let $S=\{\beta\in {\mathbb F}:\,
\beta^{2}\not \in -{\mathbb F}^{2}\}$, we now have a relation
on $S$ given by 
       $$\tilde{\beta}\sim \beta\mbox{\  if and only if\ }
        \tilde{\beta}/\beta\in G(\beta)$$ 
From Proposition 2.8, we know that this is an equivalence relation. We can also see this directly. First notice that 
   $$(a^{2}+b^{2}\beta)^{2}=(a^{2}+(b/c)^{2}\beta c^{2})^{2}$$
for all $c\in {\mathbb F}^{*}$. Hence $G(\beta)=G(\beta c^{2})$ for all 
$c\in {\mathbb F}^{*}$. In particular we have that $G(\tilde{\beta})=
G(\beta)$ if $\tilde{\beta}\sim \beta$. Let us now see that $\sim$ is
an equivalence relation. Firstly it is refexive as $\beta/\beta=1\in G(\beta)$,
the latter being a group. To see that $\sim$ is symmetric, notice that 
$G(\tilde{\beta})=G(\beta)$ is a group and thus 
$\tilde{\beta}/\beta\in G(\beta)$ if and only if $\beta/\tilde{\beta}\in
G(\tilde{\beta})$. Finally to see that $\sim$ is transitive, let $\alpha,\beta,
\delta\in S$ such that $\alpha\sim\beta$ and $\beta\sim\delta$. Then
$G(\alpha)=G(\beta)=G(\delta)$ and $\beta/\alpha,\delta/\beta \in
G(\alpha)$ implies that their product $\delta/\alpha\in G(\alpha)$. \\ \\
Let us now move to the case when the characteristic of ${\mathbb F}$ is
$2$. We first see that the algebras here split naturally into two classes.
\begin{lemm} Let $L$ be a nilpotent SAA of dimension $10$ with an isotropic center of dimension $4$ where $L^3=Z(L)$ and $L$ is of type C over a field ${\mathbb F}$ of 
characteristic $2$. Then $L$ cannot have both a presentation of
the form ${\mathcal P}_{10}^{(4,4)}(0,\beta)$ and ${\mathcal P}_{10}^{(4,4)}(\alpha,\gamma)$ where $\alpha\not =0$ and $t^2+\beta, t^2+\alpha t+\gamma$ are irreducible polynomials in ${\mathbb F}[t]$.
\end{lemm}
{\bf Proof}\ \ We argue by contradiction and suppose we have
an algebra satisfying both types of presentations. By Proposition 2.6
we then have $0=(ad+bc)\alpha=(ad-bc)\alpha=\alpha$ that contradicts
our assumption that $\alpha\not =0$. $\Box$ \\ \\
For the algebras where $\alpha =0$, the same analysis works as for algebras
with $\mbox{char\,}{\mathbb F}\not =2$ and we get the same result as in Proposition 2.8. This leaves us with algebras where  $\mbox{char\,}{\mathbb F}=2$
and where $\alpha\not =0$. Notice that Proposition 2.6 tells us here that
${\mathcal P}_{10}^{(4,4)}(\alpha,\beta)$ and ${\mathcal P}_{10}^{(4,4)}(\tilde{\alpha},\tilde{\beta})$
describe the same algebra if and only if 
\begin{eqnarray*}
     \tilde{\alpha} & = & \frac{\alpha}{d^{2}+ cd \alpha+c^{2}\beta} \\
    \tilde{\beta} & = & \frac{b^{2}+ab\alpha +a^{2}\beta}{d^{2}+cd\alpha +
c^{2}\beta}
\end{eqnarray*}
for some $a,b,c,d\in {\mathbb F}$ where $ad+bc=1$. \\ \\ 
We don't take the analysis further but end by considering an example, the
finite fields of characteristic $2$. Let ${\mathbb F}$ be the finite
field of order $2^{n}$. As a first step we show that we can always in that case,
choose our presentation such that $\beta =1$. To see this take
first some arbitrary $\alpha$ and $\beta$ such that $L$ satisfies the
presentation ${\mathcal P}_{10}^{(4,4)}(\alpha,\beta)$. The groups of units ${\mathbb F}^{*}$
is here a cyclic group of odd order $2^{n}-1$ and thus $({\mathbb F}^{*})^{2}
={\mathbb F}^{*}$.
Now pick $b\in {\mathbb F}^{*}$ such that $b^{2}=\beta/\alpha$ and let $a=0$, $c=1/b$
and $d=b$. Then $ad+bc=1$ and 
   $$\tilde{\beta}=\frac{b^{2}}{b^{2}+\alpha +\beta/b^{2}}=1.$$ 
Thus we can assume from now on that $\beta=1$. 
Now let $b\in {\mathbb F}$ be arbitrary and let $a=b+1,c=b$ and $d=b+1$. Then 
$ad+bc=1$, $\tilde{\beta}=1$ and $\tilde{\alpha}=\frac{\alpha}{b(b+1)\alpha +1}$. The number of such values for $\tilde{\alpha}$ is $2^{n-1}$ that gives us
all possible values such that $t^{2}+\tilde{\alpha}t+1$ is irreducible (easy
to count that the number of reducible polynomials of the form $t^{2}+
ut+1$ is $2^{n-1}$). We thus conclude that there is only one algebra of
type C in this case. 
\section{Algebras with an isotropic centre of dimension $2$}
In this section we will be assuming that $Z(L)$ is isotropic of dimension $2$. 
Notice that if $L={\mathbb F}u+{\mathbb F}v+L^{2}$, then $L^{2}={\mathbb F}uv+
L^{3}$. It follows that $L^{2}=Z(L)^{\perp}$ is of dimension $8$ and that 
$L^{3}$ is of dimension $7$. We can then pick our standard basis such that 
\begin{eqnarray*}
     Z(L) & = & {\mathbb F}x_{5}+{\mathbb F}x_{4}, \\
   L^{2} & = & {\mathbb F}x_{5}+\cdots +{\mathbb F}x_{1}+{\mathbb F}y_{1}+
               {\mathbb F}y_{2}+{\mathbb F}y_{3}, \\
  L^{3} & = & {\mathbb F}x_{5}+\cdots + {\mathbb F}x_{1}+{\mathbb F}y_{1}+
              {\mathbb F}y_{2}.
\end{eqnarray*}
Furthermore $L^{3}={\mathbb F}uvu+{\mathbb F}uvv+L^{4}$ and thus 
 $\mbox{dim\,}L^{4}\in \{5,6\}$. Let $k$ be the nilpotence class of $L$. We know that
the maximal class is $7$ and as $\mbox{dim\,}L^{k}\not =1$ and $L^{k}\leq
Z(L)$, we must have that $L^{k}=Z(L)$. Moreover, we know that $\mbox{dim\,}L^{s}\not =2$ for $1\leq s\leq 4$ and thus $5\leq k\leq 7$. If $L^{5}=Z(L)$ then $\mbox{dim\,}Z_{2}(L)-\mbox{dim\,}Z(L)=\mbox{dim\,}L^{2}-\mbox{dim\,}L^{3}=1$ and we get
the contradiction that $L^{4}=Z_{2}(L)$ is of dimension $3$. Thus $6\leq k\leq 7$. We will deal with the two cases separately.
\subsection{The algebras of class $6$}
As the class is $6$, it follows that $(L^{4},L^{4})=(L^{7},L)=0$ and thus
$L^{4}$ is isotropic. We have seen that the dimension of $L^{4}$ is at least
$5$ and thus $\mbox{dim\,}L^{4}=5$. We can thus now furthermore choose
our standard basis such that \\
$$\begin{array}{ll}
\begin{tabular}{c|c|c|c} 
\cline{2-3} 
  \mbox{} & $x_{5}$ & $y_{5}$ & \mbox{} \\
 \mbox{}$L^{6}=Z(L)$ & $x_{4}$ & $y_{4}$ &  \mbox{} \\
\cline{2-3}
\mbox{}$L^{5}=Z_{2}(L)$  & $x_{3}$ & $y_{3}$ & $L^{2}=Z_{5}(L)$ \\ 
\cline{2-3}
 \mbox{} & $x_{2}$ & $y_{2}$ & $L^{3}=Z_{4}(L)$ \\
  \mbox{}$L^{4}=Z_{3}(L)$ & $x_{1}$ & $y_{1}$ & \mbox{} \\
\cline{2-3}
\end{tabular} &  
$$\begin{array}{l}
 Z(L)=L^{6}={\mathbb F}x_{5}+{\mathbb F}x_{4} \\ 
 Z_{2}(L)=L^{5}={\mathbb F}x_{5}+{\mathbb F}x_{4}+{\mathbb F}x_{3} \\
 Z_{3}(L)=L^{4}={\mathbb F}x_{5}+\cdots + {\mathbb F}x_{1} \\
Z_{4}(L)=L^{3}=L^{4}+{\mathbb F}y_{1}+{\mathbb F}y_{2} \\
Z_{5}(L)=L^{2}=L^{4}+{\mathbb F}y_{1}+{\mathbb F}y_{2}+{\mathbb F}y_{3}. 
\end{array}$$
\end{array}$$
As $L^{4}Z_{4}(L)=0$ we must have 
   $$x_{1}y_{2}=0.$$
It then follows that $L^{4}L^{3}=0$ and then 
   $$L^{3}L^{3}={\mathbb F}y_{1}y_{2}.$$ 
Notice that $L^{3}L^{3}\not =0$ since this would imply that $(L^{6},L)=
(L^{3},L^{4})=(L^{3}L^{3},L)=0$ and we would get the contradiction that
the class of $L$ is at most $5$. Next let us see that $x_{1}y_{3}$ and
$x_{2}y_{3}$ are linearly independent. To see this we argue by contradiction
and suppose that $(ax_{1}+bx_{2})y_{3}=0$ for some $a,b\in {\mathbb F}$ where
not both $a,b$ are zero. But this would imply that  
$(ax_{1}+bx_{2})L\in Z(L)$ and we would thus get the contradiction that 
$ax_{1}+bx_{2}\in Z_{2}(L)=L^{5}$. Thus 
      $$L^{4}L^{2}={\mathbb F}x_{1}y_{3}+{\mathbb F}x_{2}y_{3}=Z(L).$$
Notice that $L^{3}L^{3}={\mathbb F}y_{1}y_{2}$ is a one-dimensional characteristic subspace of $Z_{2}(L)$. We consider two possibilities: $L^{3}L^{3}\leq Z(L)$
and $L^{3}L^{3}\not\leq Z(L)$.  
\subsubsection{Algebras where $L^{3}L^{3}\leq Z(L)$}
We pick our standard basis such that 
\begin{equation}
      L^{3}L^{3}={\mathbb F}y_{1}y_{2}={\mathbb F}x_{5}.
\end{equation}
We have seen above that $Z(L)={\mathbb F}x_{2}y_{3}+{\mathcal F}x_{1}y_{3}=
L^{4}L^{2}$. In order to clarify the structure of $L$ we introduce
that characteristic subspace 
     $$W=\{x\in L^{4}:\, xL^{2}\leq L^{3}L^{3}\}.$$
Notice that $W$ is the kernel of the surjective linear map $L^{4}\rightarrow Z(L)/L^{3}L^{3},\,x\mapsto xy_{3}+L^{3}L^{3}$ and thus $W$ is of codimension $1$ in $L^{4}$. Also
$L^{5}<W$. We can thus pick our standard basis such that 
    $$W={\mathbb F}x_{5}+{\mathbb F}x_{4}+{\mathbb F}x_{3}+{\mathbb F}x_{2}.$$
From this one sees that  we have a chain of characteristic ideals
of $L$ \\ \\
\begin{tabular}{c|c|c|c} 
\cline{2-3} 
  $L^{3}L^{3}$ & $x_{5}$ & $y_{5}$ & \mbox{} \\
\cline{2-3}
 \mbox{}$L^{6}=Z(L)$ & $x_{4}$ & $y_{4}$ &  $(L^{3}L^{3})^{\perp}$ \\
\cline{2-3}
\mbox{}$L^{5}=Z_{2}(L)$  & $x_{3}$ & $y_{3}$ & $L^{2}=Z_{5}(L)$ \\ 
\cline{2-3}
 \mbox{}$W$ & $x_{2}$ & $y_{2}$ & $L^{3}=Z_{4}(L)$ \\
\cline{2-3}
  \mbox{}$L^{4}=Z_{3}(L)$ & $x_{1}$ & $y_{1}$ & $W^{\perp}$ \\
\cline{2-3}
\end{tabular} 
\mbox{}\\ \\ \\
Notice that ${\mathbb F}x_{2}y_{3}={\mathbb F}x_{5}$. We continue considering 
characteristic subspaces. Let 
        $$S=\{x\in L^{3}:\,x\cdot L^{2}\leq L^{3}L^{3}\}.$$
Notice that $L^{3}L^{2}=Z(L)$ and that $S$ is the kernel of the surjective
linear map $L^{3}\rightarrow Z(L)/L^{3}L^{3},\,x\mapsto x\cdot y_{3}+L^{3}L^{3}$
and is thus of codimension $1$ in $L^{3}$. Notice also that $x_{1}\not\in S$
whereas $W\leq S$. 
It follows that we can pick our standard basis such that
 $$ S={\mathbb F}x_{5}+{\mathbb F}x_{4}+{\mathbb F}x_{3}+{\mathbb F}x_{2}+
{\mathbb F}y_{1}+{\mathbb F}y_{2}.$$
In particular we have $y_{1}y_{3},y_{2}y_{3}\in L^{3}L^{3}$. Notice that 
    $$S^{\perp}=L^{5}+{\mathbb F}y_{1}$$
and that $L^{2}S^{\perp}=L^{2}y_{1}={\mathbb F}y_{1}y_{2}+{\mathbb F}y_{1}y_{3}=L^{3}L^{3}$. Let 
       $$T=\{x\in L^{2}:\,xS^{\perp}=0\}.$$
Then $T$ is the kernel of the surjective linear map $L^{2}\rightarrow 
L^{3}L^{3},\,x\mapsto y_{1}x$. Notice that $W^{\perp}\leq T$ but that $y_{2}\not\in T$. We can then pick our standard basis such that
     $$T=W^{\perp}+{\mathbb F}y_{3}.$$
In particular $y_{1}y_{3}=0$. We now have a characteristic isotropic 
subspace $T^{\perp}={\mathbb F}x_{5}+{\mathbb F}x_{4}+{\mathbb F}x_{2}$ where
$T^{\perp}\cdot (L^{3}L^{3})^{\perp}=x_{2}(L^{3}L^{3})^{\perp}=L^{3}L^{3}$. We now
let 
     $$R=\{x\in (L^{3}L^{3})^{\perp}:\,xT^{\perp}=0\}.$$
This is the kernel of the surjective linear map $(L^{3}L^{3})^{\perp}\rightarrow
L^{3}L^{3},\,x\mapsto x_{2}x$ that contains $L^{3}$. We now  refine our standard basis such that 
     $$R=L^{3}+{\mathbb F}y_{4}$$
and we have in particular $x_{2}y_{4}=0$. Let us summarize. For every standard
basis that respects the list of characteristic subspaces above, we have
\begin{eqnarray*}
    x_{2}y_{4} & = & 0 \\  
  \mbox{}{\mathbb F}x_{2}y_{3} & = & {\mathbb F}x_{5}  \\
     x_{1}y_{2} & = & 0 \\ 
{\mathbb F}x_{5}+{\mathbb F}x_{1}y_{3} & = & {\mathbb F}x_{5}+
{\mathbb F}x_{4} \\  
 {\mathbb F}x_{5}+{\mathbb F}x_{1}y_{4} & = & {\mathbb F}x_{5}+{\mathbb F}x_{3}  \\
  {\mathbb F}y_{1}y_{2} & = & {\mathbb F}x_{5} \\ 
  y_{1}y_{3} & = & 0 \\ 
 y_{2}y_{3} & \in & {\mathbb F}x_{5}
\end{eqnarray*}
It is not difficult to see that we can furthermore refine our basis such that 
\begin{equation}
  x_{2}y_{3}=x_{5},\ x_{2}y_{4}=0,\ x_{1}y_{2}=0,\ x_{1}y_{3}=x_{4},\ x_{1}y_{4}=-x_{3},\ 
   y_{1}y_{2}=x_{5}, y_{1}y_{3}=0.
\end{equation}  
This deals with all triple values apart from
$$\begin{array}{lll}
(y_{1}y_{4},y_{5})=b, & (y_{2}y_{4},y_{5})=e, & (x_{3}y_{4},y_{5})=r, \\
(y_{3}y_{4},y_{5})=f, & (y_{2}y_{3},y_{5})=c, & \mbox{} 
\end{array}$$
Replacing $x_1, y_1, x_2, y_2, x_4, y_4, x_5, y_5, y_3$ by $r^2 x_1, (1/r)^2 y_1,  (1/r) x_2 + (b/r) x_4 , r y_2 - cr x_2 - e x_3 - bcr x_4,  r^2 x_4, (1/r)^2 y_4 - (b/r^2) y_2, (1/r) x_5 , ry_5 , y_3 - (e/r) x_2 - (1/r)(f+bc) x_3 - (b/r)e x_4$ and we can assume that $b=f=e=c=0$ and $r=1$.
We thus arrive at a unique presentation. 
\begin{prop} There is a unique nilpotent SAA $L$ of dimension $10$ with an
isotropic center of dimension $2$ that has the further properties
that $L$ is nilpotent of class $6$ and $L^{3}L^{3}\leq Z(L)$. This algebra
can be given by the presentation
 $${\mathcal P}_{10}^{(2,1)}:\ \ (x_{3}y_{4},y_{5})=1,\ (x_{2}y_{3},y_{5})=1,\ 
      (x_{1}y_{3},y_{4})=1,\ (y_{1}y_{2},y_{5})=1.$$
\end{prop}
One readily verifies that the algebra with the presentation above has the properties stated. 
\subsubsection{Algebras where $L^{3}L^{3}\nleq Z(L)$}
Recall that $x_{1}y_{2}=0$. As we had observed before, ${\mathbb F}y_{1}y_{2}=L^{3}L^{3}\leq Z_{2}(L)=Z(L)+{\mathbb F}x_{3}$. We had also seen that ${\mathbb F}x_{1}y_{3}+
{\mathbb F}x_{2}y_{3}=Z(L)$. We can now pick our standard basis such that
\begin{equation}
   x_{1}y_{2}=0,\  y_{1}y_{2}=x_{3},\ x_{1}y_{3}=x_{4},\ x_{2}y_{3}=x_{5}.
\end{equation}
This leaves us with the following list of triple values to determine.
$$\begin{array}{llll}
  (x_{1}y_{4},y_{5})=a,  & (y_{2}y_{3},y_{4})=c, & (y_{1}y_{3},y_{4})=f, & 
(y_{1}y_{4},y_{5})=h, \\
\mbox{}  (x_{2}y_{4},y_{5})=b,  & (y_{2}y_{3},y_{5})=d, & (y_{1}y_{3},y_{5})=g, &
(y_{3}y_{4},y_{5})=k, \\
 \mbox{}  (x_{3}y_{4},y_{5})=r, & (y_{2}y_{4},y_{5})=e, & \mbox{} &
\end{array}$$
where $r\, \not =\, 0$ as $x_{3} \,\not \in \, Z(L)$. We show that we can further refine the basis so that $a=b=c=d=e=f=g=h=k=0$. Let $\alpha = ac-e+bd,$ $\beta = c-g$ and replace $x_1, y_1, x_2, y_2,$ $y_3, y_4, y_5$ by  
$x_1 - (a/r) x_3, y_1 - c x_2 - ((h-bc)/r) x_3, x_2 - (b/r)x_3 - f x_4 + \beta x_5, 
y_2 - c x_1 - d x_2 + (\alpha/r) x_3,$
$y_3 - (h/r) x_1 - (e/r) x_2 - (k/r) x_3 + (a/r)y_1 + (b/r)y_2$, 
$y_4 - cf x_1 - df x_2 + (\alpha f/r) x_3 + f y_2$, 
$y_5 + c \beta x_1 + d \beta x_2 - (\beta \alpha/r) x_3 - \beta y_2$.
We thus arrive at a family of algebras given by
the presentation ${\mathcal P}_{10}^{(2,2)}$ given in the next Proposition.
\begin{prop}
Let $L$ be a nilpotent SAA of dimension $10$ with an isotropic center of 
dimension $2$ with the further properties that $L$ is of nilpotence class $6$
and $L^{3}L^{3}\not\leq Z(L)$. Then $L$ has a presentation of the form
     $${\mathcal P}_{10}^{(2,2)}(r):\ \ (x_{3}y_{4},y_{5})=r,\ 
      (x_{2}y_{3},y_{5})=1,\ (x_{1}y_{3}y_{4})=1,\ (y_{1}y_{2},y_{3})=1,$$
where $r\not =0$. Furthermore the presentations ${\mathcal P}_{10}^{(2,2)}(r)$
and ${\mathcal P}_{10}^{(2,2)}(s)$ describe the same algebra if and only
if $s/r\in ({\mathbb F}^{*})^{4}$. Conversely any algebra with such a presentation has the properties stated.
\end{prop}
{\bf Proof}\ \ We have already seen that all such algebras have a presentation
of the form ${\mathcal P}_{10}^{(2,2)}(r)$ for some $0\not =r\in {\mathbb F}$.
Straightforward calculations show that conversely any algebra with such a presentation has the properties stated in the Proposition. It remains to prove
the isomorphism property. To see that it is sufficient, suppose that
we have an algebra $L$ with presentation $\ {\mathcal P}_{10}^{(2,2)}\,(r)\ $ with
respect to some given standard basis. Let $\, s\, $ be any element in ${\mathbb F}^{*}$ such that $s/r \, = \, b^{4}\, \in \, ({\mathbb F}^{*})^{4}$. Replace the basis for $L$
with a new standard basis $\tilde{x}_{1},\cdots, \tilde{y}_{5}$ where 
$\tilde{x}_{1}=x_{1}$, $\tilde{y}_{1}=y_{1}$, $\tilde{x}_{2}=(1/b)x_{2},
\tilde{y}_{2}=by_{2}$, 
$\tilde{x}_{3}=bx_{3}$, 
$\tilde{y}_{3}=(1/b)y_{3}$,
$\tilde{x}_{4}=(1/b)x_{4}$,
$\tilde{y}_{4}=by_{4}$, 
$\tilde{x}_{5}=
(1/b^{2})x_{5}$
and 
$\tilde{y}_{5}=b^{2}y_{5}$. Direct calculations show that
$L$ has the presentation ${\mathcal P}_{10}^{(2,2)}(s)$ with respect to this
new basis.
It remains to see that the condition is necessary. Consider again an algebra
$L$ with presentation ${\mathcal P}_{10}^{(2,2)}(r)$ and suppose that
$L$ has also a presentation ${\mathcal P}_{10}^{(2,2)}(s)$ with respect to
some other standard basis $\tilde{x}_{1},\ldots ,\tilde{y}_{5}$. We want
to show that $s/r\in ({\mathbb F}^{*})^{4}$. We know that 
$L={\mathbb F}\tilde{y}_{5}+{\mathbb F}\tilde{y}_{4}+L^{2}=
{\mathbb F}y_{5}+{\mathbb F}y_{4}+L^{2}$. Thus 
\begin{eqnarray*}
    \tilde{y}_{4} & = & ay_{4}+by_{5} + u_{4} \\
    \tilde{y}_{5} & = & cy_{4}+dy_{5}+u_{5},
\end{eqnarray*}
for some $u_{4},u_{5}\in L^{2}$ and $a,b,c,d\in {\mathbb F}$ where 
$ad-bc\not =0$. We know that $L^{2}L^{2}\leq L^{4}$ and thus 
\begin{eqnarray*}
    \tilde{y}_{5}\tilde{y}_{4}\tilde{y_{4}} & = & 
                       (cy_{4}+dy_{5})(ay_{4}+by_{5})(ay_{4}+by_{5})+w \\
    \tilde{y}_{4}\tilde{y}_{5}\tilde{y}_{5} & = & 
              (ay_{4}+by_{5})(cy_{4}+dy_{5})(cy_{4}+dy_{5})+z,
\end{eqnarray*}
where $w,z\in L^{4}$. We use the fact that $L^{7}=0$ and $L^{3}L^{3}\leq L^{5}$ in the 
following calculation. We have 
      $$s^{3}=(\tilde{y}_{4}\tilde{y}_{5}\tilde{y}_{5}
  ((\tilde{y}_{4}\tilde{y}_{5}),\tilde{y}_{5}\tilde{y}_{4}\tilde{y}_{4})=
      (ad-bc)^{4}r^{3}.$$
Hence $s/r\in ({\mathbb F}^{*})^{4}$. $\Box$ \\ \\
{\bf Remarks}.(1)  It thus depends on the field ${\mathbb F}$, how many algebras there are of this type. When $({\mathbb F}^{*})^{4}={\mathbb F}^{*}$ there is just
one algebra. This includes the case when ${\mathbb F}$ is algebraically closed
or a finite field of characteristic $2$. \\ \\
(2) Let ${\mathbb F}$ be a finite field of
order $p^{n}$ where $p$ is an odd prime. If $p\equiv 1\,(\mbox{mod}\,4)$ then
there are $4$ algebras and if $p\equiv -1\,(\mbox{mod}\,4)$ then there are
$4$ algebras when $n$ is even and $2$ algebras when $n$ is odd. \\ \\
(3) For ${\mathbb F}={\mathbb R}$ there are two algebras, one for $r<0$ and
one for $r>0$. For ${\mathbb F}={\mathbb Q}$ there are infinitely many
algebras. 
\subsection{The algebras of class $7$}
Here we are dealing with algebras of maximal class and thus we can make use
of the general theory concerning these. In particular we know that we can choose our standard basis such that
\begin{eqnarray*}
L^{7}=Z(L) & = & {\mathbb F}x_{5}+{\mathbb F}x_{4}, \\
\mbox{}L^{6}=Z_{2}(L) & = & {\mathbb F}x_{5}+{\mathbb F}x_{4}+{\mathbb F}x_{3}, \\
\mbox{}L^{5}=Z_{3}(L) & = & {\mathbb F}x_{5}+{\mathbb F}x_{4}+
     {\mathbb F}x_{3}+{\mathbb F}x_{2}.
\end{eqnarray*}
We also know that $L^{4}=Z_{4}(L)=(L^{5})^{\perp}$, $L^{3}=Z_{5}(L)=(L^{6})^{\perp}$
and $L^{2}=Z_{6}(L)=(L^{7})^{\perp}$. Furthermore we know that can also get
characteristic ideals of dimensions $1,5$ and $9$ in the following way.
 \\ \\
Firstly, we know from the general theory that $x_{3}y_{4},x_{2}y_{3}\not =0$. As a result $L^{5}L^{2}={\mathbb F}x_{2}y_{3}\not =0$. This gives us a characteristic ideal of dimension $1$ and then $(L^{5}L^{2})^{\perp}$ is a characteristic
ideal of dimension $9$. \\ \\
We now turn to the description of a characteristic ideal of dimension $5$.  From
the general theory we  also know that $x_{1}y_{2},y_{2}y_{2}$ are linearly independent. Thus  
$L^{4}L^{3}={\mathbb F}x_{1}y_{2}+{\mathbb F}y_{1}y_{2}$ is a $2$-dimensional
characteristic subspace of $L^{6}$. Let $I_{1}=L^{5}L^{2}, I_{2}=L^{7}$
and $I_{3}=L^{6}$. Let $k$ be smallest such that $I_{k}\cap L^{4}L^{3}\not =\{0\}$. Then 
        $$U=\{x\in L^{4}:\,xL^{3}\leq I_{k}\}$$
is a characteristic ideal of dimension $5$. We can thus further refine
our basis such that we have the following situation. \\ \\
\begin{tabular}{c|c|c|c} 
\cline{2-3} 
  $L^{5}L^{2}$ & $x_{5}$ & $y_{5}$ & \mbox{} \\
\cline{2-3}
 \mbox{}$L^{7}=Z(L)$ & $x_{4}$ & $y_{4}$ &  $(L^{5}L^{2})^{\perp}$ \\
\cline{2-3}
\mbox{}$L^{6}=Z_{2}(L)$  & $x_{3}$ & $y_{3}$ & $L^{2}=Z_{6}(L)$ \\ 
\cline{2-3}
 \mbox{}$L^{5}=Z_{3}(L)$ & $x_{2}$ & $y_{2}$ & $L^{3}=Z_{5}(L)$ \\
\cline{2-3}
  \mbox{}$U$ & $x_{1}$ & $y_{1}$ & $L^{4}=Z_{4}(L)$ \\
\cline{2-3}
\end{tabular} 
\mbox{}\\ \\ \\
There are now few separate cases to consider according to whether $L^{4}L^{3}=L^{7}$
or $L^{4}L^{3}\not =L^{7}$ and whether or not $L^{4}L^{3}\cap L^{5}L^{2}\not =\{0\}$. 
\subsubsection{Algebras where $L^{4}L^{3}=L^{7}$}
 In this case we have
     $${\mathbb F}x_{2}y_{3}={\mathbb F}x_{1}y_{2}={\mathbb F}x_{5},\ 
     {\mathbb F}x_{5}+{\mathbb F}y_{1}y_{2}={\mathbb F}x_{5}+{\mathbb F}x_{4}.$$
Now consider the characteristic subspace $UL^{2}={\mathbb F}x_{5}+{\mathbb F}x_{1}y_{3}$. There are again two subcases to consider as either $UL^{2}$ has dimension $1$ or $2$. \\ \\
{\bf I. Algebras where $UL^{2}$ is $1$-dimensional} \\ \\
In this case we have that $x_{1}y_{3}\in {\mathbb F}x_{5}$ and it follows that
$L^{4}L^{2}=UL^{2}+y_{1}L^{2}={\mathbb F}x_{5}+y_{1}L^{2}=Z(L)$. Consider
the characteristic subspace
    $$V=\{x\in L^{2}:\,L^{4}x\leq L^{5}L^{2}\}.$$
Then $V$ is of codimension $1$ in $L^{2}$ and $L^{4}\leq V$. Also $y_{2}\not \in
V$. One sees readily that we can refine our choice of basis further such that
     $$V=L^{4}+{\mathbb F}y_{3}.$$
In particular $y_{1}y_{3}\leq {\mathbb F}x_{5}$. Next consider the characteristic
subspace 
    $$W=\{x\in L^{4}:\,xV=0\}.$$
We have that $L^{6}\leq W$ and that $x_{2}\not\in W$. Also $W$ is the kernel
of the surjective linear map $L^{4}\rightarrow L^{5}L^{2},\,x\mapsto xy_{3}$
and thus of codimension $1$ in $L^{4}$. We can now pick our basis further such
that 
       $$W=L^{6}+{\mathbb F}x_{1}+{\mathbb F}y_{1}.$$
It is not difficult to see that such a choice is compatible to what we have
done so far. Notice that it follows that $y_{1}y_{3}=x_{1}y_{3}=0$. Next one notices that $L^{3}V=Z(L)$ and that $L^{3}U\leq L^{5}L^{2}$. Let 
     $$Z=\{x\in V:\,L^{3}x\leq L^{5}L^{2}\}.$$
Then $Z$ is of codimension $1$ in $V$ and $y_{1}\not\in Z$. We can now further refine the basis such that 
      $$Z=U+{\mathbb F}y_{3}.$$
The reader can convince himself that this is compatible to our choice
so far. In particular $y_{2}y_{3}\in {\mathbb F}x_{5}$. 
%
%
%
%
%
%
Replacing $y_{2}$ by a suitable $y_{2}-ux_{2}$, we can furthermore assume that
$y_{2}y_{3}=0$. 
With this choice of basis we thus have $x_{1}y_{3}=y_{1}y_{3}=y_{2}y_{3}=0$
as well as ${\mathbb F}x_{1}y_{2}={\mathbb F}x_{2}y_{3}={\mathbb F}x_{5}$
and ${\mathbb F}x_{5}+{\mathbb F}x_{4}={\mathbb F}x_{5}+{\mathbb F}y_{1}y_{2}$.
It is not difficult to see that we can further refine our basis such that
\begin{equation}
  x_{1}y_{3}=y_{1}y_{3}=y_{2}y_{3}=0,\ x_{1}y_{2}=x_{2}y_{3}=x_{5},\ y_{1}y_{2}=x_{4}.
\end{equation}
We are then only left with the following triple values
$$\begin{array}{lll}
   (x_{1}y_{4},y_{5})=a, &  (y_{1}y_{4},y_{5})=c, & (y_{3}y_{4},y_{5})=e, \\ 
   (x_{2}y_{4},y_{5})=b, &  (y_{2}y_{4},y_{5})=d, & (x_{3}y_{4},y_{5})=r,
\end{array}$$
where $r\not =0$. 
Replacing $x_{5},x_{4},x_3, x_2, x_1, y_{1},y_2, y_{3},y_{4},y_{5}$ by
$(1/r)^2 x_5,
r^4 x_4,
r x_3 + br x_4,
(1/r) x_2 + (a/r) x_4 - (c/r) x_5,
(1/r)^3 x_1,
r^3 y_1 + d r^3 x_4,
ry_2,
(1/r)y_3 - (e/r^2) x_3 - (e/r^2)b x_4,
(1/r)^4 y_4 -(b/r^4) y_3 - (a/r^4) y_2 + (d/r^4) x_1,
r^2 y_5 + c r^2 y_2
$, gives
us a new standard basis where we can assume that $a=b=c=d=e=0$ and $r=1$. The reader can
check that (12) is not affected by these changes. We thus
arrive at unique presentation for $L$.

\begin{prop}
There is a unique nilpotent SAA $L$ of dimension $10$ that has isotropic 
center of dimension $2$ with the further properties that the class is
$7$, $L^{4}L^{3}=L^{7}$ and $\mbox{dim\,}UL^{2}=1$. This algebra can
be given by the presentation 
       $${\mathcal P}_{10}^{(2,3)}:\ (x_{3}y_{4},y_{5})=1,\ 
    (x_{2}y_{3},y_{5})=1,\ (x_{1}y_{2},y_{5})=1,\ (y_{1}y_{2},y_{4})=1.$$
\end{prop}
Direct calculations show that the algebra with this presentation has the properties stated. \\ \\
{\bf II. Algebras where $UL^{2}$ is $2$-dimensional} \\ \\
In this case ${\mathbb F}x_{2}y_{3}={\mathbb F}x_{1}y_{2}={\mathbb F}x_{5}$
and ${\mathbb F}y_{1}y_{2}+{\mathbb F}x_{5}={\mathbb F}x_{1}y_{3}+
{\mathbb F}x_{5}={\mathbb F}x_{4}+{\mathbb F}x_{5}$. It is not difficult
to see that we can choose our standard basis such that 
\begin{equation}
  x_{1}y_{2}=x_{2}y_{3}=x_{5},\ y_{1}y_{2}=x_{4}.
\end{equation}
Now $y_{1}y_{3}=ax_{5}+bx_{4}$ for some $a,b\in {\mathbb F}$. Replacing 
$x_{2}$, $y_{1}$, $y_{2}$, $y_{3}$ by $x_{2}+bx_{3}$, $y_{1}-ax_{2}$,
$y_{2}-ax_{1}$, $y_{3}-by_{2}+abx_{1}$ gives 
\begin{equation}
     y_{1}y_{3}=0,
\end{equation}
and the changes do not affect (13). Next consider $y_{2}y_{3}=ax_{4}+bx_{5}$
and replace $x_{1}$, $y_{2}$, $y_{3}$ by $x_{1}-ax_{3}$, $y_{2}-bx_{2}$, 
$y_{3}+ay_{1}$. These changes imply that we can assume furthermore that
\begin{equation}
    y_{2}y_{3}=0.
\end{equation}
Now consider $x_{1}y_{3}=ax_{5}+bx_{4}$ (where $b\not =0$ by our assumptions). Replacing $x_{1},y_{2}$ by $x_{1}-ax_{2},y_{2}+ay_{1}$ we can assume that $a=0$. Then
replace $x_{1},\ldots ,y_{5}$ by $(1/b)x_{1}$, $(1/b^{2})x_{2}$, $(1/b^{3})x_{3}$,
$b^{3}x_{4}$, $bx_{5}$, $by_{1}$, $b^{2}y_{2}$, $b^{3}y_{3}$, $(1/b^{3})y_{4}$,
$(1/b)y_{5}$ we can assume that $b=1$. Thus
\begin{equation}
x_{1}y_{3}=x_{4}.
\end{equation}
This leaves us with the following triples.
$$\begin{array}{lll}
   (x_{1}y_{4},y_{5})=a, &  (y_{1}y_{4},y_{5})=c, & (y_{3}y_{4},y_{5})=e, \\ 
   (x_{2}y_{4},y_{5})=b, &  (y_{2}y_{4},y_{5})=d, & (x_{3}y_{4},y_{5})=r,
\end{array}$$
where $r \neq 0$. We refine our basis further such that $a=b=c=d=e=0$. Replacing $y_5, y_1, x_3, y_3,$
$ x_2, y_4$ by 
$y_5 + c y_2 , 
y_1 + d x_4 , 
x_3 + b x_4,
y_3 - (e/r) x_3 - b (e/r) x_4 , 
x_2 + a x_4 - c x_5 , 
y_4 - a y_2 - b y_3 + d x_1$.
Thus $L$ has a presentation of the form ${\mathcal P}_{10}^{(2,4)}(r)$ as described
in the next proposition. 
\begin{prop}
Let $L$ be a nilpotent SAA of dimension $10$ with an isotropic center of
dimension $2$ that is of class $7$ and has the further properties
that $L^{4}L^{3}=L^{7}$ and $\mbox{dim\,}UL^{2}=2$. This algebra can
be given by a presentation of the form
   $${\mathcal P}_{10}^{(2,4)}(r):\ (x_{3}y_{4},y_{5})=r,\ 
       (x_{2}y_{3},y_{5})=1,\ (x_{1}y_{2},y_{5})=1,\ (x_{1}y_{3},y_{4})=1,\ (y_{1}y_{2},y_{4})=1,$$
where $r\not =0$. Furthermore two such presentations ${\mathcal P}_{10}^{(2,4)}(r)$ and ${\mathcal P}_{10}^{(2,4)}(s)$ describe the same algebra if and only
if $s/r\in ({\mathbb F}^{*})^{11}$. Conversely any algebra with such a presentation has the properties stated.
\end{prop}
{\bf Proof}\ \ \ We have already seen that any such algebra has such a presentation. Direct calculations show that an algebra with a presentation 
$\, \, {\mathcal P}_{10}^{(2,4)}(r)\, $ has the properties stated. We turn to the isomorphism property. To see that the condition is sufficient, suppose we have 
an algebra $L$ with a presentation ${\mathcal P}_{10}^{(2,4)}(r)$ with 
respect to some standard basis $x_{1},y_{1},\ldots,x_{5},y_{5}$. Suppose that 
$s/r=a^{11}$ for some $a\in {\mathbb F}^{*}$. Consider a new standard
basis $\tilde{x}_{1}=ax_{1}$, 
$\tilde{y}_{1}=(1/a)y_{1}$,
$\tilde{x}_{2}=a^{3}x_{2}$, 
$\tilde{y}_{2}=(1/a^{3})y_{2}$, 
$\tilde{x}_{3}=
a^{5}x_{3}$, 
$\tilde{y}_{3}=(1/a^{5})y_{3}$, 
$\tilde{x}_{4}=(1/a^{4})=x_{4}$,
$\tilde{y}_{4}=a^{4}y_{4}$, 
$\tilde{x}_{5}=(1/a^{2})x_{5}$, 
$\tilde{y}_{5}=
a^{2}y_{5}$. 
Calculations show that $L$ has then presentation ${\mathcal P}_{10}^{(2,4)}(s)$ with respect to the new basis.
It only remains now to see that the conditions is also necessary. Consider an algebra $L$ with presentation ${\mathcal P}_{10}^{(2,4)} (r)$ with respect to
some standard basis $x_{1},y_{1},\ldots ,x_{5},y_{5}$. Take some arbitrary
new standard basis $\tilde{x}_{1},\tilde{y}_{1}, \ldots ,\tilde{x}_{5},
\tilde{y}_{5}$ such that $L$ satisfies the presentation
 ${\mathcal P}_{10}^{(2,4)} (s)$ with respect to the new basis. Using the fact that
we have an ascending chain of characteristic ideals we know that 
$$\begin{array}{l}
  \tilde{y}_{1}=ay_{1}+\beta_{11}x_{1}+\cdots + \beta_{15}x_{5}, \\
  \tilde{y}_{2}=by_{2}+\alpha_{21}y_{1}+\beta_{21}x_{1}+\cdots +
\beta_{25}x_{5}, \\
  \tilde{y}_{3}=cy_{3}+\alpha_{32}y_{2}+\alpha_{31}y_{1}+\beta_{31}x_{1}+\cdots +
\beta_{35}x_{5}, \\
  \tilde{y}_{4}=dy_{4}+\alpha_{43}y_{3}+\alpha_{42}y_{2}+\alpha_{41}y_{1}+
    \beta_{41}x_{1}+\cdots +\beta_{45}x_{5}, \\
\tilde{y}_{5}=ey_{5}+\alpha_{54}y_{4}+\cdots +\alpha_{51}y_{1}+
\beta_{51}x_{1}+\cdots +\beta_{55}x_{5}, \\
 \tilde{x}_{1}=(1/a)x_{1}+\gamma_{12}x_{2}+\cdots +\gamma_{15}x_{5}, \\
 \tilde{x}_{2}=(1/b)x_{2}+\gamma_{23}x_{3}+\gamma_{24}x_{4}+\gamma_{25}x_{5}, \\
 \tilde{x}_{3}=(1/c)x_{3}+\gamma_{34}x_{4}+\gamma_{35}x_{5}, \\
\tilde{x}_{4}=(1/d)x_{4}+\gamma_{45}x_{5}, \\
\tilde{x}_{5}=(1/e)x_{5},
\end{array}$$
for some $\alpha_{ij},\beta_{ij},\gamma_{ij},a,b,c,d,e\in {\mathbb F}$ where
$a,b,c,d,e\not =0$. 
Direct calculations show that 
\begin{eqnarray*}
   1=(\tilde{x}_{1}\tilde{y}_{2},\tilde{y}_{5})=be/a & \Rightarrow &
           e=a/b \\
\mbox{} 1=(\tilde{y}_{1}\tilde{y}_{2},\tilde{y}_{4})=abd & \Rightarrow &
        d=1/(ab) \\ 
\mbox{} 1=(\tilde{x}_{1}\tilde{y}_{3},\tilde{y}_{4})=cd/a & \Rightarrow &
      c=a^{2}b \\
\mbox{} 1=(\tilde{x}_{2}\tilde{y}_{3},\tilde{y}_{5})=ce/b & \Rightarrow &
      b=a^{3}.
\end{eqnarray*}
Thus $b=a^{3}$, $c=a^{5}$, $d=1/a^{4}$, $e=1/a^{2}$ and it follows that
   $$s=(\tilde{x}_{3}\tilde{y}_{4},\tilde{y}_{5})=(de/c)r=(1/a)^{11}r.$$
Hence $s/r\in ({\mathbb F}^{*})^{11}$. $\Box$ \\ \\
{\bf Remarks}. It follows that if $({\mathbb F}^{*})^{11}={\mathbb F}^{*}$ then
there is only one algebra of this type. This includes any algebraically
closed field and ${\mathbb R}$. If ${\mathbb F}$ is a finite field
of order $p^{n}$, then the number of algebras is either $11$ or $1$ according to
whether $11$ divides $p^{n}-1$ or not. Notice also that there are infinitely
many algebras over ${\mathbb Q}$. 
\subsubsection{Algebras where $L^{4}L^{3}\neq L^{7}$ and $L^{5}L^{2}\leq L^{4}L^{3}$}
Here we can pick our basis such that 
     $$L^{4}L^{3}={\mathbb F}x_{5}+{\mathbb F}x_{3}.$$
Notice also that as before $U=\{x\in L^{4}:xL^{3}\leq L^{5}L^{2}\}$ and
thus again $x_{1}y_{2}\in {\mathbb F}x_{5}$.
Notice also that 
     $$(L^{4}L^{3})^{\perp}={\mathbb F}x_{5}+\cdots + {\mathbb F}x_{1}+
{\mathbb F}y_{1}+{\mathbb F}y_{2}+{\mathbb F}y_{4}.$$
Then $L^{5}(L^{4}L^{3})^{\perp}=({\mathbb F}x_{3}+{\mathbb F}x_{2})y_{4}=
L^{5}L^{2}$. Consider the characteristic subspace 
    $$V=\{x\in L^{5}:\, x(L^{4}L^{3})^{\perp}=0\}.$$
Here $x_{3}y_{4}\not =0$ and thus $V$ is the kernel of a surjective linear
map $L^{5}\rightarrow L^{5}L^{2},\ x\mapsto xy_{4}$ and has codimension $1$
in $L^{5}$. We pick our standard basis such that 
     $$V={\mathbb F}x_{5}+{\mathbb F}x_{4}+{\mathbb F}x_{2}.$$
In particular 
\begin{equation}
x_{2}y_{4}=0.
\end{equation}
Here we have again $UL^{2}={\mathbb F}x_{5}+{\mathbb F}x_{1}y_{3}$ and thus
either the dimension of $UL^{2}$ is $1$ or $2$. We consider these cases separately. \\ \\
{\bf I. Algebras where $UL^{2}$ is $1$-dimensional} \\ \\
Notice that 
$$V^{\perp}={\mathbb F}x_{5}+ \cdots +{\mathbb F}x_{1}+{\mathbb F}y_{1}+
{\mathbb F}y_{3}.$$
and that $UV^{\perp}={\mathbb F}x_{5}=L^{5}L^{2}$. Let 
       $$W=\{x\in U:\,xV^{\perp}=0\}.$$
Here $x_{2}y_{3}\not =0$ and $W$ is the kernel of the surjective linear map
$U\rightarrow L^{5}L^{2},\ x\mapsto xy_{3}$. We choose our standard basis
further such that 
     $$W={\mathbb F}x_{5}+{\mathbb F}x_{4}+{\mathbb F}x_{3}+{\mathbb F}x_{1}.$$
In particular 
\begin{equation}
    x_{1}y_{3}=0.
\end{equation}
Next look at $L^{4}V^{\perp}={\mathbb F}x_{5}+{\mathbb F}y_{1}y_{3}$. Notice
that $y_{1}y_{3}\in V$ and that $(y_{1}y_{3},y_{2})\not = 0$ (as ${\mathbb F}y_{1}y_{2}+ {\mathbb F}x_{5}={\mathbb F}x_{3}+{\mathbb F}x_{5}$). We choose our basis
further such that 
       $$L^{4}V^{\perp}={\mathbb F}x_{5}+{\mathbb F}x_{2}.$$
In particular
\begin{equation}
    (y_{1}y_{3},y_{4})=0.
\end{equation}
Now consider the characteristic subspace
       $$T=L^{4}V^{\perp}+L^{4}L^{3}={\mathbb F}x_{5}+{\mathbb F}x_{3}+
{\mathbb F}x_{2}.$$
Notice that $T^{\perp}=L^{4}+{\mathbb F}y_{4}$ and $WT^{\perp}={\mathbb F}x_{3}y_{4}+{\mathbb F}x_{1}y_{4}=L^{5}L^{2}$. Let 
     $$R=\{x\in W:\, xT^{\perp}=0\}.$$
We have $x_{3}y_{4}\not =0$ and $R$ is the kernel of the surjective linear
map $W\rightarrow L^{5}L^{2},\ x\mapsto xy_{4}$. We now refine our basis
further such that
       $$R={\mathbb F}x_{5}+{\mathbb F}x_{4}+{\mathbb F}x_{1}.$$
In particular $x_{1}y_{4}=0$. We have thus got a basis where $x_{1}y_{3}=x_{1}y_{4}=x_{2}y_{4}=0$ and where $(y_{1}y_{3},y_{4})=0$. It is not difficult that
we can furthermore assume that
$$\begin{array}{l}
   x_{3}y_{4}=x_{5},\ 
   x_{2}y_{4}=0,\ x_{1}y_{2}=x_{5}, \\
  x_{1}y_{3}=0,\ x_{1}y_{4}=0,\ y_{1}y_{2}=x_{3},\ (y_{1}y_{3},y_{4})=0.
\end{array}$$
We still need to consider the following triples.
$$\begin{array}{lll}
   (y_{3}y_{4},y_{5})=a, &  (y_{1}y_{4},y_{5})=c, & (y_{2}y_{4},y_{5})=f, \\ 
   (y_{1}y_{3},y_{5})=b, &  (y_{2}y_{3},y_{4})=d, & (y_{2}y_{3},y_{5})=e, \\
   (x_{2}y_{3},y_{5})=r, &        &
\end{array}$$
where $r \neq 0$. Replacing $x_1, x_2, y_1, y_2, y_3, y_4,y_5$ by 
$x_1 + d x_4$,
$x_2 - b x_5$,
$y_1 - c x_3$,
$y_2 - (1/r)(c+e) x_2 - f x_3 + (b/r)(c+e) x_5$,
$y_3 - c x_1 - f x_2 - (a+bd) x_3 - cd x_4 + bf x_5$,
$y_4 -d y_1$,
$y_5 + b y_2$ we can assume that $a=b=c=d=e=f =0$.
It follows that $L$ has a presentation of the form ${\mathcal P}_{10}^{(2,5)}(r)$
as in the following proposition.
\begin{prop}
Let $L$ be a nilpotent SAA of dimension $10$ with an isotropic center of 
dimension $2$ that is of class $7$ and has the further properties that 
$L^{4}L^{3}\not =L^{7}$, $L^{5}L^{2}\leq L^{4}L^{3}$ and $UL^{2}$ is $1$-dimensional. The algebra can be given by a presentation of the form
$${\mathcal P}_{10}^{(2,5)}(r):\ (x_{2}y_{3},y_{5})=r,\ (x_{3}y_{4},y_{5})=1,\
 (x_{1}y_{2},y_{5})=1,\ (y_{1}y_{2},y_{3})=1,$$
where $r\not =0$. Furthermore two such presentations 
${\mathcal P}_{10}^{(2,5)}(r)$ and ${\mathcal P}_{10}^{(2,5)}(s)$ describe the
same algebra if and only if $s/r\in ({\mathbb F}^{*})^{3}$. Conversely any algebra with such a presentation has the properties stated.
\end{prop}
{\bf Proof}\ \ We have already seen that any such algebra has a presentation
of this form. Conversely, direct calculations show that any algebra
with a presentation ${\mathcal P}_{10}^{(2,5)}(r)$ satisfies the
properties stated. We turn to the isomorphism property. To see that condition
is sufficient, suppose we have an algebra $L$ with presentation 
${\mathcal P}_{10}^{(2,5)}(r)$ with respect to some standard basis
$x_{1},y_{1},\ldots ,x_{5},y_{5}$. Suppose that $s/r=a^{3}$ for some
$a\in {\mathbb F}^{*}$. Consider a new standard basis $\tilde{x}_{1}=x_{1}$,
$\tilde{y}_{1}=y_{1}$, $\tilde{x}_{2}=ax_{2}$, $\tilde{y}_{2}=(1/a)y_{2}$, 
$\tilde{x}_{3}=(1/a)x_{3}$, $\tilde{y}_{3}=ay_{3}$, $\tilde{x}_{4}=x_{4}$,
$\tilde{y}_{4}=y_{4}$, $\tilde{x}_{5}=(1/a)x_{5}$, $\tilde{y}_{5}=ay_{5}$. Calculations show that $L$ has then the presentation ${\mathcal P}_{10}^{(2,5)}(s)$
with respect to the new basis. 
It only remains to see that the condition is also necessary. Consider an
algebra $L$ with presentation ${\mathcal P}_{10}^{(2,5)}(r)$ with respect
to some standard basis $x_{1},y_{1},\ldots ,$
$x_{5}$
$,y_{5}$. Take some arbitrary new
standard basis $\tilde{x}_{1},\tilde{y}_{1}, \ldots ,\tilde{x}_{5},
\tilde{y}_{5}$ such that $L$ also satisfies presentation 
${\mathcal P}_{10}^{(2,5)}(s)$ with respect to the new basis. Using the fact
that we have an ascending chain of characteristic ideals as well
as the fact that $L^{4}L^{3}={\mathbb F}x_{5}+{\mathbb F}x_{3}$, $L^{4}V^{\perp}=
{\mathbb F}x_{5}+{\mathbb F}x_{2}$, $R={\mathbb F}x_{5}+{\mathbb F}x_{4}+
{\mathbb F}x_{1}$, $W^{\perp}={\mathbb F}y_{2}+U$, $(R+L^{4}V^{\perp})^{\perp}=
{\mathbb F}y_{3}+U$ and $T^{\perp}={\mathbb F}y_{4}+{\mathbb F}y_{1}+U$
 are characteristic subcases, we know that 
$$\begin{array}{l}
  \tilde{y}_{1}=(1/a)y_{1}+\beta_{11}x_{1}+\cdots + \beta_{15}x_{5}, \\
  \tilde{y}_{2}=(1/b)y_{2}+\beta_{21}x_{1}+\cdots +
\beta_{25}x_{5}, \\
  \tilde{y}_{3}=(1/c)y_{3}+\beta_{31}x_{1}+\cdots +
\beta_{35}x_{5}, \\
  \tilde{y}_{4}=(1/d)y_{4}+\alpha_{41}y_{1}+
    \beta_{41}x_{1}+\cdots +\beta_{45}x_{5}, \\
\tilde{y}_{5}=(1/e)y_{5}+\alpha_{54}y_{4}+\cdots +\alpha_{51}y_{1}+
\beta_{51}x_{1}+\cdots +\beta_{55}x_{5}, \\
 \tilde{x}_{1}=ax_{1}+\gamma_{14}x_{4}+\gamma_{15}x_{5}, \\
 \tilde{x}_{2}=bx_{2}+\gamma_{25}x_{5}, \\
 \tilde{x}_{3}=cx_{3}+\gamma_{35}x_{5}, \\
\tilde{x}_{4}=dx_{4}+\gamma_{45}x_{5}, \\
\tilde{x}_{5}=ex_{5},
\end{array}$$
for some $a,b,c,d,e,\alpha_{ij},\beta_{ij},\gamma_{ij}$ where $a,b,c,d,e\not =0$.
It follows that 
$$\begin{array}{l}
  1=(\tilde{x}_{3}\tilde{y}_{4},\tilde{y}_{5})=c/(de) \\
\mbox{}  1=(\tilde{x}_{1}\tilde{y}_{2},\tilde{y}_{5})=a/(be) \\
\mbox{}  1=(\tilde{y}_{1}\tilde{y}_{2},\tilde{y}_{3})=1/(abc).
\end{array}$$
This gives  $c=1/(ab),\ e=a/b,\ d=1/a^{2}$ and then 
    $$s=(\tilde{x}_{2}\tilde{y}_{3},\tilde{y}_{5})=br/(ce)=b^{3}r.$$ 
This finishes the proof. $\Box$ \\ \\
{\bf Remarks}. Again we just get one algebra if $({\mathbb F}^{*})^{3}=
{\mathbb F}^{*}$. This includes all fields that are algebraically
closed as well as ${\mathbb R}$. For a finite field of order $p^{n}$ there
are $3$ algebras if $3|p^{n}-1$ but otherwise one. For ${\mathbb Q}$
there are infinitely many algebras. \\ \\
{\bf I. Algebras where $UL^{2}$ is $2$-dimensional} \\ \\
Recall that $UL^{2}={\mathbb F}x_{5}+{\mathbb F}x_{1}y_{3}$, 
$L^{4}L^{3}={\mathbb F}x_{5}+{\mathbb F}x_{3}$ and $x_{2}y_{4}=0$. It is not difficult to see that one can further refine the basis such that
\begin{equation}
   x_{3}y_{4}=\alpha x_{5},\ x_{2}y_{3}=rx_{5},\ x_{2}y_{4}=0,\ x_{1}y_{2}=x_{5},\ x_{1}y_{3}=
x_{4},\ y_{1}y_{2}=x_{3},
\end{equation}
where $\alpha,r \neq 0$. Replacing $x_{1},y_{1},\ldots ,x_{5},y_{5}$ by $\alpha x_{1}$, 
$(1/\alpha)y_{1}$, $(1/\alpha^{3})x_{2}$, $\alpha^{3}y_{2}$,
$\alpha^{2}x_{3}$, $(1/\alpha^{2})y_{3}$, $(1/\alpha)x_{4}$,
$\alpha y_{4}$, $\alpha^{4}x_{5}$, $(1/\alpha^{4})y_{5}$, implies
that we can furthermore assume that $\alpha =1$. We have also the following
triples to sort out.
$$\begin{array}{lll}
   (y_{1}y_{3},y_{4})=a, &  (y_{2}y_{3},y_{4})=d, & (x_{1}y_{4},y_{5})=g, \\ 
   (y_{1}y_{3},y_{5})=b, &  (y_{2}y_{3},y_{5})=e, & (y_{3}y_{4},y_{5})=h. \\
   (y_{1}y_{4},y_{5})=c, &  (y_{2}y_{4},y_{5})=f.  &
\end{array}$$
Let  $\gamma = c+e$, $\beta =ae-h-bd$ and replace $x_2, x_1, y_1, y_2, y_3,y_4, y_5$ by
$x_2 - a x_4 - b x_5, 
x_1 - (a+g) x_3 + d x_4 + \gamma  x_5, 
y_1 - c x_3, 
y_2 -f x_3,
y_3 - c x_1 - f x_2 + \beta x_3 + (af-cd) x_4 + bf x_5 + (a+g) y_1,
y_4 - d y_1 + a y_2,
y_5 + c \gamma  x_3 - \gamma  y_1 + b y_2 $.
Then we can assume that $a=b=c=d=e=f=g=h=0$ and $L$ has a presentation ${\mathcal P}_{10}^{(2,6)}(r)$
as in the following proposition.
\begin{prop}
Let $L$ be a nilpotent SAA of dimension $10$ with an isotropic center of 
dimension $2$ that is of class $7$ and has the further properties that 
$L^{4}L^{3}\not =L^{7}$, $L^{5}L^{2}\leq L^{4}L^{3}$ and $UL^{2}$ is $2$-dimensional. The algebra can be given by a presentation of the form
$${\mathcal P}_{10}^{(2,6)} (r):\ (x_{2}y_{3},y_{5})=r,\ (x_{3}y_{4},y_{5})=1,\
 (x_{1}y_{2},y_{5})=1,\ (x_{1}y_{3},y_{4})=1,\ (y_{1}y_{2},y_{3})=1$$
where $r\not =0$. Furthermore two such presentations 
${\mathcal P}_{10}^{(2,6)}(r)$ and ${\mathcal P}_{10}^{(2,6)}(s)$ describe the
same algebra if and only if $s/r\in ({\mathbb F}^{*})^{12}$. Conversely any algebra with such a presentation has the properties stated.
\end{prop}
{\bf Proof}\ \ We have already seen that any such algebra has a presentation
of this form. Conversely, direct calculations show that any algebra
with a presentation ${\mathcal P}_{10}^{(2,6)}(r)$ satisfies the
properties stated. We turn to the isomorphism property. To see that condition
is sufficient, suppose we have an algebra $L$ with presentation 
${\mathcal P}_{10}^{(2,6)}(r)$ with respect to some standard basis
$x_{1},y_{1},\ldots ,x_{5},y_{5}$. Suppose that $s/r=a^{12}$ for some
$a\in {\mathbb F}^{*}$. Consider a new standard basis $\tilde{x}_{1}=(1/a)x_{1}$,
$\tilde{y}_{1}=ay_{1}$, $\tilde{x}_{2}=a^{4}x_{2}$, $\tilde{y}_{2}=(1/a^{4})y_{2}$, 
$\tilde{x}_{3}=(1/a^{3})x_{3}$, $\tilde{y}_{3}=a^{3}y_{3}$, 
$\tilde{x}_{4}=a^{2}x_{4}$,
$\tilde{y}_{4}=(1/a^{2})y_{4}$, $\tilde{x}_{5}=(1/a^{5})x_{5}$, 
$\tilde{y}_{5}=a^{5}y_{5}$. 
Calculations show that $L$ has then the presentation $\, {\mathcal P}_{10}^{(2,6)} (s)\, $
with respect to the new basis. 
It only remains to see that the condition is also necessary. Consider an
algebra $L$ with presentation ${\mathcal P}_{10}^{(2,6)} (r)$ with respect
to some standard basis $x_{1},y_{1},\ldots ,$
$x_{5},y_{5}$. Take some arbitrary new
standard basis $\tilde{x}_{1},\tilde{y}_{1}, \ldots $
$,\tilde{x}_{5}, \tilde{y}_{5}$ such that $L$ also satisfies presentation 
${\mathcal P}_{10}^{(2,6)}(s)$ with respect to the new basis. Using the fact
that we have an ascending chain of characteristic ideals as well
as the fact that $L^{4}L^{3}={\mathbb F}x_{5}+{\mathbb F}x_{3}$, $V=
{\mathbb F}x_{5}+{\mathbb F}x_{4}+{\mathbb F}x_{2}$,
 are characteristic subspaces, we know that 
$$\begin{array}{l}
  \tilde{y}_{1}=(1/a)y_{1}+\beta_{11}x_{1}+\cdots + \beta_{15}x_{5}, \\
  \tilde{y}_{2}=(1/b)y_{2}+\alpha_{21}y_{1}+\beta_{21}x_{1}+\cdots +
\beta_{25}x_{5}, \\
  \tilde{y}_{3}=(1/c)y_{3}+\alpha_{31}y_{1}+\beta_{31}x_{1}+\cdots +
\beta_{35}x_{5}, \\
  \tilde{y}_{4}=(1/d)y_{4}+\alpha_{42}y_{2}+\alpha_{41}y_{1}+
    \beta_{41}x_{1}+\cdots +\beta_{45}x_{5}, \\
\tilde{y}_{5}=(1/e)y_{5}+\alpha_{54}y_{4}+\cdots +\alpha_{51}y_{1}+\beta_{51}x_{1}+\cdots +\beta_{55}x_{5}, \\
 \tilde{x}_{1}=ax_{1}+\gamma_{12}x_{2}+\ldots + \gamma_{15}x_{5}, \\
 \tilde{x}_{2}=bx_{2}+\gamma_{24}x_{4}+\gamma_{25}x_{5}, \\
 \tilde{x}_{3}=cx_{3}+\gamma_{35}x_{5}, \\
\tilde{x}_{4}=dx_{4}+\gamma_{45}x_{5}, \\
\tilde{x}_{5}=ex_{5},
\end{array}$$
for some $a,b,c,d,e,\alpha_{ij},\beta_{ij},\gamma_{ij}$ where $a,b,c,d,e\not =0$.
It follows that 
$$\begin{array}{l}
  1=(\tilde{x}_{3}\tilde{y}_{4},\tilde{y}_{5})=c/(de) \\
\mbox{}  1=(\tilde{x}_{1}\tilde{y}_{2},\tilde{y}_{5})=a/(be) \\
\mbox{}  1=(\tilde{y}_{1}\tilde{y}_{2},\tilde{y}_{3})=1/(abc) \\
\mbox{} 1=(\tilde{x}_{1}\tilde{y}_{3},\tilde{y}_{4})=a/(cd).
\end{array}$$
This gives  $b=(1/a^{4}),\ c=a^{3},\ d=(1/a^{2}),\ e=a^{5}$ and then 
    $$s=(\tilde{x}_{2}\tilde{y}_{3},\tilde{y}_{5})=br/(ce)=(1/a^{12})r.$$ 
This finishes the proof. $\Box$ \\ \\
{\bf Remarks}. The number of algebras depends thus again on the underlying
field. In particular there is one algebra over the field ${\mathbb C}$, 
two algebras over ${\mathbb R}$ and infinitely many over ${\mathbb Q}$.
When ${\mathbb F}$ is a finite field of order $p^{n}$ then the number can
be $12$, $6$, $3$, $4$, $2$ or $1$ depending on what the value of $p^{n}$ is modulo $12$.
\subsubsection{Algebras where $L^{4}L^{3}\neq L^{7}$ and $L^{5}L^{2}\nleq L^{4}L^{3}$}
First we pick our standard basis such that 
   $$L^{4}L^{3}\cap L^{7}={\mathbb F}x_{4},\ L^{4}L^{3}={\mathbb F}x_{4}+
{\mathbb F}x_{3}.$$
Thus in particular ${\mathbb F}x_{1}y_{2}={\mathbb F}x_{4}$. From this one
sees that $U(L^{5}L^{2})^{\perp}=L^{7}+{\mathbb F}x_{1}y_{4}$ where
$(x_{1}y_{4},y_{2})=-(x_{1}y_{2},y_{4})\not =0$. We further refine our basis
such that
    $$U\cdot (L^{5}L^{2})^{\perp}=L^{7}+{\mathbb F}x_{2}={\mathbb F}x_{5}
+{\mathbb F}x_{4}+{\mathbb F}x_{2}.$$
Then notice that $(U(L^{5}L^{2})^{\perp})L=L^{5}L^{2}+{\mathbb F}x_{2}y_{5}$,
where $(x_{2}y_{5},y_{3})=-(x_{2}y_{3},y_{5})\not =0$. We refine our basis
further such that 
    $$(U(L^{5}L^{2})^{\perp})L=L^{5}L^{2}+{\mathbb F}x_{3}={\mathbb F}x_{5}+{\mathbb F}x_{3}.$$
Notice that in particular $x_{2}x_{5}\in {\mathbb F}x_{5}+{\mathbb F}x_{3}$ and thus 
$(x_{2}x_{4},x_{5})=-(x_{2}x_{5},x_{4})=0$. We have as well $(x_{2}y_{4},y_{3})=
-(x_{2}y_{3},y_{4})=0$ and thus
\begin{equation}
   x_{2}y_{4}=0.
\end{equation}
We also have 
     $${\mathbb F}x_{3}=({\mathbb F}x_{4}+{\mathbb F}x_{3})\cap ({\mathbb F}x_{5}+
		     {\mathbb F}x_{3})=(L^{4}L^{3})\cap ((U(L^{5}L^{2})^{\perp})L).$$
Next consider the characteristic subspace
     $$V=\{x\in L^{4}:\, xL^{3}\leq (L^{4}L^{3})\cap ((U(L^{5}L^{2})^{\perp})L)\}.$$
Notice that $x_{1}y_{2}\not =0$ and thus $V$ is the kernel of a surjective linear map $L^{4}\rightarrow 
L^{4}L^{3}/{\mathbb F}x_{3},\ x\mapsto xy_{2}+{\mathbb F}x_{3}$ and thus of codimension $1$ in $L^{4}$.
Notice also that $L^{5}\leq V$. We refine our basis further such that
        $$V=L^{5}+{\mathbb F}y_{1}={\mathbb F}x_{5}+\ldots +{\mathbb F}x_{2}+{\mathbb F}y_{1}.$$
It follows in particular that 
\begin{equation}
     {\mathbb F}y_{1}y_{2}={\mathbb F}x_{3}
\end{equation}
Notice that 
     $$(U(L^{5}L^{2})^{\perp})^{\perp}={\mathbb F}x_{5}+\cdots +{\mathbb F}x_{1}+{\mathbb F}y_{1}+
		     {\mathbb F}y_{3}$$
and then $U\cdot (U(L^{5}L^{2})^{\perp})^{\perp}=L^{5}L^{2}+{\mathbb F}x_{1}y_{3}$. Now 
$(x_{1}y_{3},y_{2})=-(x_{1}y_{2},y_{3})=0$ but also $x_{1}y_{4}\in U(L^{5}L^{2})^{\perp}={\mathbb F}x_{5}+
{\mathbb F}x_{4}+{\mathbb F}x_{2}$ and thus $(x_{1}y_{3},y_{4})=-(x_{1}y_{4},y_{3})=0$. It follows
that $x_{1}y_{3}\in {\mathbb F}x_{5}=L^{5}L^{2}$. It follows that $U\cdot (U(L^{5}L^{2})^{\perp})^{\perp}=L^{5}L^{2}$. Consider the characteristic subspace  
     $$R=\{x\in U:\, x(U(L^{5}L^{2})^{\perp})^{\perp}=0\}.$$
We have that $x_{2}y_{3}\not =0$ and thus $R$ is the kernel of the surjective linear map $U\rightarrow
L^{5}L^{2},\, x\mapsto xy_{3}$. Thus $R$ is of codimension $1$ in $U$ and contains $L^{6}$. Now choose
our basis further such that 
     $$R=L^{6}+{\mathbb F}x_{1}={\mathbb F}x_{5}+{\mathbb F}x_{4}+{\mathbb F}x_{3}+{\mathbb F}x_{1}.$$
In particular 
\begin{equation}
    x_{1}y_{3}=0.
\end{equation}
Next consider $L^{4}\cdot (U(L^{5}L^{2})^{\perp})^{\perp}=L^{5}L^{2}+{\mathbb F}y_{1}y_{3}$. Notice that
$(y_{1}y_{3},y_{2})=-(y_{1}y_{2},y_{3})\not =0$. We can refine our basis further such that 
     $$L^{4}\cdot (U(L^{5}L^{2})^{\perp})^{\perp}={\mathbb F}x_{5}+{\mathbb F}x_{2}.$$ 
In particular
\begin{equation}
      (y_{1}y_{3},y_{4})=0.
\end{equation}
It is not difficult to see that we can now choose our basis such that
$$\begin{array}{llll}
   x_{3}y_{4}=x_{5}, & x_{2}y_{3}=rx_{5}, & x_{2}y_{4}=0,  & \\
   x_{1}y_{2}=x_{4}, & x_{1}y_{3}=0, & y_{1}y_{2}=x_{3}, & (y_{1}y_{3},y_{4})=0,
\end{array}$$
where $r\neq 0$. Replacing $x_{1},y_{1},\ldots ,x_{5},y_{5}$ by $(1/r)x_{1}$, $ry_{1}$, $rx_{2}$, $(1/r)y_{2}$, $x_{3}$, $y_{3}$,
$(1/r^{2})x_{4}$, $r^{2}y_{4}$, $r^{2}x_{5}$, $(1/r^{2})y_{5}$, we see that we can furthermore assume that $r=1$.
We are now only left with the triples
 $$\begin{array}{llll}
   (x_{1}y_{4},y_{5})=a, &   (y_{1}y_{4},y_{5})=d, &  (y_{2}y_{4},y_{5})=g,  &(y_{3}y_{4},y_{5})=h. \\ 
   (y_{1}y_{3},y_{5})=c, &  (y_{2}y_{3},y_{4})=e,& (y_{2}y_{3},y_{5})=f, &  
   \end{array}$$
We now show that we can further refine the basis such that all these values are zero. 
By replacing $x_2, x_1, y_1, y_2, y_3, y_4, y_5$ by
$x_2 - c x_5, 
x_1 + (c-a) x_3 + (e+d) x_4 + f x_4, 
y_1 - d x_3, 
y_2 - g x_3, 
y_3 - d x_1 - g x_2 - (h+ce)x_3 + (cg - df) x_5 + (a-c) y_1,
y_4 + d (d+e) x_3 - (e+d) y_1 ,
y_5 - f y_1 + c y_2
$, we have thus arrived at a unique presentation. 
\begin{prop}
There is a unique nilpotent SAA of dimension $10$ with an isotropic center of 
dimension $2$ that is of class $7$ and has the further properties that 
$L^{4}L^{3}\not =L^{7}$, $L^{5}L^{2}\not\leq L^{4}L^{3}$. The algebra can be given by the presentation 
$${\mathcal P}_{10}^{(2,7)}:\ (x_{2}y_{3},y_{5})=1,\ (x_{3}y_{4},y_{5})=1,\
 (x_{1}y_{2},y_{4})=1,\ (y_{1}y_{2},y_{3})=1$$
\end{prop}
{\bf Proof}\ \ We have already seen that any such algebra must have a presentation of this form and conversely direct calculations show that the algebra with the given presentation 
satisfies all the properties stated. $\Box$
%

\bibliographystyle{plainurl}
\bibliography{biblio}

\end{document}